\documentclass[11pt]
{article}

\usepackage[margin=20mm]{geometry}
\usepackage{amsmath}%
\usepackage{amsthm}%
\usepackage{amsfonts}%
\usepackage{amssymb}%
\usepackage{amscd}
\usepackage{graphicx}
\usepackage{xcolor}
\usepackage{bbm}
\usepackage[all,cmtip]{xy}
\usepackage{wrapfig}
\usepackage{enumerate}
\usepackage{subcaption}

\usepackage{color}
\definecolor{grau}{rgb}{0.2,0.2,0.2}
\usepackage[colorlinks, linkcolor=grau, citecolor=grau, urlcolor=grau, breaklinks]{hyperref}
\usepackage{breakurl}
\graphicspath{{Figures/}}

\theoremstyle{plain}
\newtheorem{theorem}{Theorem}[section]

\newtheorem{lemma}[theorem]{Lemma}
\newtheorem{corollary}[theorem]{Corollary}
\newtheorem{conjecture}[theorem]{Conjecture}

\theoremstyle{definition}
\newtheorem{definition}[theorem]{Definition}
\newtheorem{remark}[theorem]{Remark}
\newtheorem{remarks}[theorem]{Remarks}

\usepackage{bibspacing}

\newcommand{\blagojevic}{Blagojevi\'{c}}

\newcommand{\poincare}{{Poincar\'{e}}}
\newcommand{\vrecica}{Vre\'{c}ica}
\newcommand{\zivaljevic}{\v{Z}ivaljevi\'{c}}

\newcommand{\RR}{\mathbbm{R}} 
\newcommand{\ZZ}{\mathbbm{Z}} 
\newcommand{\CC}{\mathbbm{C}} 
\newcommand{\FF}{\mathbbm{F}} 
\newcommand{\wo}{\backslash} 
\newcommand{\To}{\longrightarrow} 
\newcommand{\toG}[1]{\longrightarrow_{#1}} 
\newcommand{\dd}{\textnormal{d}} 
\newcommand{\im}{\textnormal{im}} 

\newcommand{\homeq}{\simeq} 
\newcommand{\homeo}{\cong} 
\newcommand{\iso}{\cong} 
\newcommand{\bd}{\partial} 
\newcommand{\one}{\mathbbm{1}} 

\newcommand{\incl}{\hookrightarrow} 
\newcommand{\st}{\ |\ } 
\newcommand{\qm}[1]{``{#1}''} 

\newcommand{\eps}{\varepsilon}

\newcommand{\unorientedBordism}{\textnormal{MO}}
\newcommand{\id}{\textnormal{id}}
\newcommand{\vertices}{\textnormal{vert}}

\newcommand{\simplex}{\sigma}
\newcommand{\diag}{\Delta}
\newcommand{\fatdiag}{\diag^{\rm{fat}}}

\newcommand{\wt}[1]{\widetilde{#1}}

\newcommand{\textdef}[1]{\textnormal{\emph{#1}}}

\begin{document}

\title
{On the Square Peg Problem and Its Relatives}
\author{Benjamin Matschke%
}
\date{}

\maketitle

\begin{abstract}

Toeplitz's Square Peg Problem asks whether every continuous simple closed curve in the plane contains the four vertices of a square.
It has been proved for various classes of sufficiently smooth curves, some of which are dense, none of which are open. 
In this paper we prove it for several open classes of curves, one of which is also dense.
This can be interpreted in saying that the Square Peg Problem is solved for generic curves.
The latter class contains all previously known classes for which the Square Peg Problem has been proved in the affirmative.%
\footnote{Since the appearance of this article on the arXiv, some new results have been proved, see Section~\ref{secRecentProgress}.
In particular, Tao's class of curves~\cite{Tao17integrationApproachToToeplitzProblem} contains curves that do not lie in the classes given in the current article.}

%
%
%

We also prove results about rectangles inscribed in immersed curves. 
Finally, we show that the problem of finding a regular octahedron on metric 2-spheres has a ``topological counter-example'', that is, a certain test map with boundary condition exists.
\end{abstract}


\section{Introduction}

The Square Peg Problem was first posed by Otto Toeplitz~\cite{Toe11aufgabenDerAnalysisSitus} in 1911:

\begin{conjecture}[Square Peg Problem
]
\label{conjSquarePegProblem}
Every continuous embedding $\gamma:S^1\to\RR^2$ contains four points that are the vertices of a square.
\end{conjecture}

\begin{figure}[htb]
  \centering
  \begin{minipage}[b]{0.45\textwidth}
	\centering
	\input{Figures/SquarePegMitSpirale.pspdftex}
    \caption{Example for Conjecture~\ref{conjSquarePegProblem}.}
    \label{figSquarePegMitSpirale}
  \end{minipage}
\quad
  \begin{minipage}[b]{0.45\textwidth}
    \centering
    \input{Figures/curveHavingNoSquareFullyInside.pspdftex}
	\caption{Jordan curve without an inscribed square that lies completely inside.}
    \label{figCurveWithoutSquareFullyInside}
  \end{minipage}
\end{figure}


The name Square Peg Problem might be slightly misleading: We do not require the square to lie inside the curve, otherwise there are easy counter-examples, see Figure~\ref{figCurveWithoutSquareFullyInside}.


Toeplitz' problem has been solved affirmatively for various restricted classes of curves
such as convex curves and curves that are ``smooth enough'', by various authors:
Emch~\cite{Emc13squarePeg1, Emc15squarePeg2, Emc16squarePeg3} for piecewise analytic curves with only finitely many inflection points and other singularities where the left and right sided tangents at the finitely many non-smooth points exist, Hebbert~\cite{Hebbert14InscribedSquaresAndKinematicGeometry} for quadrilaterals,
Zindler~\cite{Zin21konvexeGebilde} and Christensen~\cite{Chr50kvadrat} for convex curves,
Jerrard~\cite{Jer61inscribedSquares} for analytic curves,
Nielsen--Wright~\cite{NielsenWright95rectanglesInscribedInSymmetricContinua} for curves that are symmetric across a line or about a point,
Makeev~\cite{Mak95quadsInscribedInClosedCurve} for star-shaped $C^2$-curves that intersect every circle in at most $4$ points (more generally he proved that any circular quadrilateral can be inscribed in such curves), 
Stromquist~\cite{Str89inscribedSquares} for locally monotone curves,
\vrecica--\zivaljevic~\cite{VrZi08fultonMacPhersonCompCyclohedraAndPolygonalPegProblem} for Stromquist's class of curves,
Pak~\cite{Pak09discreteAndPolyhedralGeometry} for piecewise linear curves,
Sagols--Mar{\'{\i}}n \cite{SaMa09inscribedSquaresDigitalPlane, SaMa11inscribedSquares} for similar discretisations,
and
Canta\-rella--Denne--McCleary~\cite{CDM11squarePeg} for curves with bounded total curvature without cusps and for $C^1$-curves. 
%
The strongest version so far was due to Stromquist \cite[Thm. 3]{Str89inscribedSquares} who established the Square Peg Problem for locally monotone Jordan curves.
Here a curve $\gamma:S^1\to\RR^2$ is called locally monotone if at every point $x\in S^1$ there exists a neighborhood $U$ of $x$ and a linear function $\alpha:\RR^2\to\RR$ such that $\alpha\circ\gamma|_U$ is strictly increasing.

All known proofs are based on the fact that ``generically'' the number of squares on a curve is odd, 
which can be measured in various topological ways. 
See Klee--Wagon~\cite{KlWa96problemsInPlaneGeomAndNumberTh}, Pak~\cite{Pak09discreteAndPolyhedralGeometry}, \vrecica--\zivaljevic~\cite{VrZi08fultonMacPhersonCompCyclohedraAndPolygonalPegProblem}, 
and \cite{Mat12surveyOnSquarePeg} for surveys.
For general embedded plane curves, the problem is still open.

\medskip

Related is the Rectangular Peg Problem, which was only recently proved by Greene and Lobb~\cite{GreeneLobb20rectangles}; compare with Section~\ref{secRecentProgress}\eqref{itProgressGreeneLobb}.

\begin{theorem}[{Rectangular Peg Problem, Greene--Lobb theorem}]
\label{thmRectangularPegProblem}
Every $C^\infty$ embedding $\gamma:S^1\to\RR^2$ contains four points that are the vertices of a rectangle with a prescribed aspect ratio $r>0$.
\end{theorem}

We state it for smooth curves only, since already this is (was) a hard problem. 
Equivalently one could state Theorem~\ref{thmRectangularPegProblem} for piecewise linear curves. 
Before the present article and the recent results listed in Section~\ref{secRecentProgress}\eqref{itProgressTao}--\eqref{itProgressGreeneLobb}, it was only known to hold in the case $r=1$, that is, for inscribed squares.
The proof in Griffiths's paper~\cite{Gri91topSquarePeg} contains unfortunately an error in the calculation of intersection numbers, see \cite{Mat08diplomaThesis} for details.
The difficulty comes from the fact that, counted with orientations, every smooth curve inscribes generically zero rectangles of a prescribed aspect ratio, see Section~\ref{secRectanglesCSTMfails}.
E.g.\ an ellipse inscribes two rectangles with opposite orientations.

\medskip

\noindent
\emph{Awards.}
The author has put 100 Euros each of these two problems; only for Conjecture~\ref{conjSquarePegProblem} the bounty is still available.

\subsection{The main results}

1.) 
In Section~\ref{secSquaresOnCurves} we essentially prove the Square Peg Problem for two new classes of curves.

\begin{wrapfigure}{r}{0.2\textwidth}
\vspace{-5pt}
\centering
\includegraphics[scale=0.5]{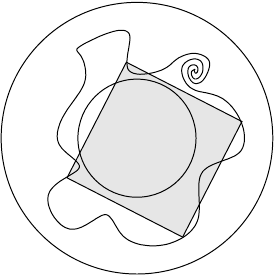}
\caption{Example \\ for Theorem~\ref{thmSquarePegInAnnulus}.}
\label{figSquarePegInAnnulus}
\vspace{-15pt}
\end{wrapfigure}

The first one is the first known \emph{open} set of continuous maps $S^1\to\RR^2$ in the compact-open topology (equivalently, in $((\RR^2)^{S^1},||.||_\infty)$) for which the Square Peg Problem holds.
It does neither require the curve to be smooth nor injective; see Section~\ref{secSquaresOnCurves} for the rather simple proof and some variations of the statement.

\begin{theorem}
\label{thmSquarePegInAnnulus}
Let $A$ denote the annulus $\{x\in\RR^2\st 1\leq ||x|| \leq 1+\sqrt{2}\}$.
Suppose that $\gamma:S^1\to A$ is a continuous closed curve in $A$ that represents the generator of $\pi_1(A)=\ZZ$.
Then $\gamma$ inscribes a square of side length at least $\sqrt{2}$.
\end{theorem}

Figure~\ref{figSquarePegInAnnulus} shows an example. 
This class does \emph{not} contain all previous known classes of curves for which the Square Peg Problem is proved, and it is not even dense, but it is the first result that bounds the size of an inscribed square from below.

2.) 
The second class is considerably more technical, but it is open, dense, and it contains all previous classes.

\begin{theorem}
\label{thmSquarePegMainTheorem}
The Square Peg Problem holds for all curves in an explicit open and dense neighbourhood of Stromquist's locally monotone curves in the space of all injective maps $S^1\to\RR^2$ with respect to the $C^0$-topology.
In this sense, any generic Jordan curve inscribes a square.
\end{theorem}

Explicit versions of this theorem are given in Corollaries~\ref{cor1MySquarePeg}, \ref{cor2MySquarePeg}, and~\ref{corSquarePegMainTheorem}.
Moreover, in Corollary~\ref{cor4MySquarePeg} we prove that curves $C^0$-close to $C^2$-embedded curves inscribe squares, which can be regarded as a version of Theorem~\ref{thmSquarePegInAnnulus} for more general shapes. 

3.) 
In Section~\ref{secRectanglesOnCurves} we prove the following first non-square special case of the Rectangular Peg Problem.

\begin{theorem}
\label{thmRectangularPegForSqrt3}
Let $\gamma:S^1\to \RR^2$ be a $C^\infty$ curve whose angular convexity is at most $60^\circ$. 
Then $\gamma$ inscribes a rectangle with aspect ratio $\sqrt{3}$.
\end{theorem}

Here we say that a smooth plane curve $\gamma$ has \emph{angular convexity} at most $\alpha$, if the signed curvature of $\gamma$ restricted to any positively oriented arc is at least $-\alpha$; here the signs are chosen in such a way that a positively oriented arc on the unit circle of length $\alpha$ has signed curvature~$+\alpha$; see Figure~\ref{figAngularConvexity}.
A smooth curve $\gamma$ is convex if and only if its angular convexity is at most~$0$.
The proof of Theorem~\ref{thmRectangularPegForSqrt3} uses a hidden symmetry that appears for~$r=\sqrt{3}$, which is a geometric piece of information.

4.)
In Section~\ref{secImmersedCurves} we deal with immersed planar curves and the parity of their inscribed squares. 
Cantarella \cite{Can08webpageOnSquarePeg} conjectured that this parity is an isotopy invariant and he stated a precise formula based on examples. 
We disprove Cantarella's conjecture and state in Theorem~\ref{thmSquaresOnImmersedCurves} how the parity can be computed from the angles at the intersection points. 
Theorem~\ref{thmRectanglesOnImmersedCurves} gives a similar formula for the parity of inscribed rectangles of a fixed aspect ratio.

\medskip


The last section, Section \ref{secCrosspolytopesOnSpheres}, treats higher-dimensional analogs.
We ask for inscribed $d$-dimensional regular crosspolytopes in metric $(d-1)$-spheres.
The problem is open for all $d\ge3$, but we use Koschorke's obstruction theory \cite{Kos81vectorFieldsAndOtherVectorBundleMorph} to derive that for $d=3$, a natural topological approach for a proof fails:
The strong test map in question exists.

%
%

%
%

Parts of this paper have been obtained in~\cite[Chap. III]{Mat08diplomaThesis} and~\cite[Chap. II]{Mat11phd}. Some of the new results have been announced in \cite{Mat09squarePegExtAbstract}.








\subsection{Recent progress} \label{secRecentProgress}

After this article appeared on the arXiv, some progress was made with respect to both Conjectures~\ref{conjSquarePegProblem} and~\ref{thmRectangularPegProblem}.
\begin{enumerate}
\item Pettersson--Tverberg--\"Osterg{\aa}rd~\cite{OPT13noteOnToeplitzSquareProblem} showed that any Jordan curve in the $12\times 12$ square grid inscribes a square whose size is at least $1/\sqrt{2}$ times the size of the largest axis parallel square that fits into the interior of the curve.

\item Van Heijst~\cite{vanHeijst14masterThesis} proved that algebraic curves of degree $d$ inscribe at most $(d^4=5d^2+4d)/4$ squares or infinitely many.
For this he makes use of Bernstein’s theorem, which states that the number of common zeros in $(\CC^*)^k$ of $k$ generic Laurent polynomials in $k$ variables with prescribed Newton polytopes equals the mixed volume of these polytopes.

\item\label{itProgressTao} Tao~\cite{Tao17integrationApproachToToeplitzProblem} proved that Jordan curves given as the union of two graphs of $(1-\eps)$-Lipschitz functions $f,g:[0,1]\to\RR$ with $f(0)=g(0)$ and $f(1)=g(1)$ inscribe squares. He makes use of an area argument, using certain conserved integrals of motion when moving squares in the plane, which appeared in a similar form already in Karasev~\cite{KaVo10MakeevsConj}.
Two special features of his approach are that he uses the Lipschitz condition as a \emph{global} condition, whereas all previous results (except for Theorem~\ref{thmSquarePegInAnnulus}) only argue locally, and that he can avoid the topological parity argument for the number of inscribed squares, which so far always creates technical problems for non-smooth curves.
If the Lipschitz constant is suitably adjusted, Tao's proof also extends to inscribed rectangles and more generally equilateral trapezoids \cite[Rem.~3.10]{Tao17integrationApproachToToeplitzProblem}.

\item The area argument by Karasev and Tao was then used by Akopyan and Avvakumov~\cite{AkopyanAvvakumov17CyclicQuads} and the author~\cite{Mat18quadrilateralsInConvexCurves} to prove Theorem~\ref{thmRectangularPegProblem} in the special case of convex curves, also for more general inscribed quadrilaterals.

\item 
\label{itProgressGreeneLobb}
Greene and Lobb \cite{GreeneLobb20rectangles} proved Theorem~\ref{thmRectangularPegProblem} (thus winning one of the two above awards).
Their proof is based on symplectic geometry and uses in particular that the Klein bottle does not admit a Lagrangian embedding into $\CC^2$.
In \cite{GreeneLobb20circularQuads} they extended this furthermore to all circular quadrilaterals, thus settling a conjecture of Makeev~\cite{Mak95quadsInscribedInClosedCurve, Mak05quadsInscribedInCurveAndVerticesOfCurve}.
\end{enumerate}

\section{Inscribed squares} \label{secSquaresOnCurves}

\subsection{Notations and a convenient parameter space of inscribed polygons} \label{secNoation}

For an element $x$ of the circle $S^1\homeo\RR/_\ZZ$ and $t\in\RR$ we define $x+t\in S^1$ as the counter-clockwise rotation of $x$ by the angle $2\pi t$ around $0$. For any space $X$, we denote by $\diag_{X^n}:=\{(x,\ldots,x)\in X^n\}$ the thin diagonal of $X^n$. Let $\simplex^n=\{(t_1,\ldots,t_{n+1})\in\RR_{\geq 0}^{n+1}\st \sum t_i=1\}$ be the standard $n$-simplex. 

The natural \textdef{parameter space} of $n$-gons is 
\[
P_n:=S^1\times \simplex^{n-1}. 
\]
It parametrises $n$-gons on $S^1$ or on some given curve $S^1\to\RR^\infty$ by their vertices in the following way
\[
\textstyle
\varphi:\ P_n\to (S^1)^n:\ (x;t_1,\ldots,t_n)\mapsto \Big(x,x+t_1,x+t_1+t_2,\ldots,x+\sum\limits_{i=1}^{n-1} t_i\Big).
\]
The so parametrised polygons are the ones that are lying counter-clockwise on $S^1$.
The map $\varphi$ is not injective, as all $(x;0,\ldots,0,1,0,\ldots,0)$ are mapped to the same point $(x,\ldots,x)$; but it is injective on $P_n\wo(S^1\times\vertices(\simplex^{n-1}))$, and on this set $\varphi$ bijects to $(S^1)^n\wo\diag_{(S^1)^n}$. 
Let $P_n^\circ=S^1\times (\simplex^{n-1})^\circ$ denote the interior of $P_n$.
The map $\varphi$ identifies $P_n^\circ$ with the set of $n$-tuples of pairwise distinct points in counter-clockwise order on $S^1$.
We define the boundary as $\bd P_n^\circ:=P_n\wo P_n^\circ$.

We let $\ZZ/n=\langle\eps\rangle$ act on $P_n$ by
\[
\eps\cdot(x;t_1,\ldots,t_n) = (x+t_1;t_2,\ldots,t_n,t_1). 
\]
This corresponds to a cyclic relabeling of the vertices of the parametrised polygon.
It makes $P_n$ clearly into a free $\ZZ/n$-space.

\begin{remark}[Other natural coordinates]
$\ZZ/n$ acts on $S^1$ via $\eps\cdot x=x+\frac{1}{n}$ and on $\simplex^{n-1}$ via $\eps\cdot(t_1,\ldots,t_n)=(t_2,\ldots,t_n,t_1)$.
Then $P_n$ is in fact as a $\ZZ/n$-space isomorphic to the product $S^1\times\simplex^n$ (with the diagonal $\ZZ/n$-action).
A particularly convenient isomorphism is the map $P_n\toG{\ZZ/4} S^1\times\simplex^n$ defined by
\[
\textstyle
(x;t_1,\ldots,t_n)\mapsto \Big(x+\sum\limits_{i=1}^n \dfrac{n+1-2i}{2n}t_i; t_1,\ldots,t_n\Big),
\]
since it sends $(x;\frac{1}{n},\ldots,\frac{1}{n})$ to $(x;\frac{1}{n},\ldots,\frac{1}{n})$.
%
%
\end{remark}

An \emph{arc} on $S^1$ from $x\in S^1$ to $y\in S^1$ will always mean the arc that goes counter-clockwise, and for $x=y$ it is degenerate to a point.
For $x,y\in S^1$, we denote by $y-x$ the length of the arc from $x$ to $y$, normalised with the factor $\frac{1}{2\pi}$.
For an $n$-tuple $(x_1,\ldots,x_n)\in\varphi(P_n)\subseteq (S^1)^n$ we write
\[
\textstyle
[x_1,\ldots,x_n]:=\Big(x_1;x_2-x_1,x_3-x_2,\ldots,x_n-x_{n-1},1-\sum\limits_{k=2}^n (x_k-x_{k-1})\Big)\in P_n.
\]
The function $[\ldots]: \varphi(P_n)\to P_n$ is right-inverse to $\varphi$, but not continuous.

\emph{Smooth} will always mean $C^\infty$.
An \emph{$\eps$-close square} is a quadrilateral whose ratios between the edges and diagonals are up to an $\eps$-error the ones of a square. 
The precise definition will not matter. 
We will use ``$\eps$-closeness'' with other polygons analogously.

\subsection{Schnirelman's proof for the smooth Square Peg Problem} \label{secShnirelmansProof}

We start with Schnirelman's proof \cite{Shn44geomPropClosedCurves}, since it is in my point of view the most beautiful one, and we will extend it in Section~\ref{secSquarePegAsimpleOpenCondition}.
The following presentation uses transversality and a bordism argument; in Schnirelman's days, these notions had not been formalised and baptised yet, but his argument works like this.

\begin{proof}
Suppose that $\gamma$ is smooth. $P_4$ parametrises quadrilaterals on $\gamma$.
Let $f:P_4\to\RR^6$ be the function that measures the four edges and the two diagonals of the quadrilaterals,
\begin{equation}
\label{eqTestMapSquarePeg}
\begin{aligned}
f:\ &P_4 &\To&\ \ \RR^4\times\RR^2 \\
&[x_1,x_2,x_3,x_4]&\longmapsto&\ (||\gamma(x_1)-\gamma(x_2)||,||\gamma(x_2)-\gamma(x_3)||,||\gamma(x_3)-\gamma(x_4)||,\\
&&&\ \ ||\gamma(x_4)-\gamma(x_1)||,||\gamma(x_1)-\gamma(x_3)||,||\gamma(x_2)-\gamma(x_4)||)
\end{aligned}
\end{equation}
Let $\Delta:=\diag_{\RR^4}\times \diag_{\RR^2}=\{(a,a,a,a,b,b)\in\RR^6\}$.
The map $f$ measures squares, since $Q:=f^{-1}(\Delta)\wo\diag_{(S^1)^4}=f^{-1}(\Delta)\cap P_4^\circ$ is the set of all squares that lie counter-clockwise on~$\gamma$.
Moreover $f$ is $\ZZ/4$-equivariant with respect to the natural $\ZZ/4$-actions.
We can deform $f$ relative to a small neighborhood of $\bd P_4^\circ$ equivariantly by a small $\eps$-homotopy to make it transversal to $\Delta$.
So $Q$ becomes a free $\ZZ/4$-manifold that parametrises $\eps$-close squares.
As the codimension of $\Delta$ in $\RR^6$ agrees with $\dim P_4=4$, $Q$ is zero-dimensional.
If we deform the curve smoothly to another curve, e.g.\ the ellipse, which can also be performed in $\RR^4$ to construct such a homotopy easily, then $Q$ changes by a $\ZZ/4$-bordism.
If the homotopy is chosen smoothly then this bordism stays away from the boundary of $P_4^\circ$, since then no curve inscribes $\eps$-close squares which have arbitrarily small edges (the angles get too close to~$\pi$).
Hence $Q$ represents a unique class $[Q]$ in the zero-dimensional unoriented bordism group $\unorientedBordism_0(P_4^\circ/_{\ZZ/4})\iso H_0(P_4^\circ/_{\ZZ/4};\ZZ/2)\iso\ZZ/2$.
If $\gamma$ is an ellipse then $f$ is already transversal to $\Delta$ and $Q$ consists of one point.
Hence $[Q]$ is the generator of $\ZZ/2$, thus $Q$ is non-empty for any smooth curve $\gamma$.
Taking a convergent subsequence of $\eps$-close squares finishes the proof.
\end{proof}

If $\gamma$ is only continuous one might try to approximate it with smooth curves and then take a convergent subsequence of the squares that are inscribed in them.
The problem is to guarantee that this subsequence does not converge to a square that degenerates to a point.
Natural candidates for which this works are continuous curves with bounded total curvature without cusps, see Cantarella, Denne \& McCleary \cite{CDM11squarePeg}.
So far, nobody managed to do this for all continuous curves.

Schnirelman's proof can be refined to get a slightly stronger result. 
\begin{corollary}[of the proof]
We may assume that $\gamma$ is parametrised counter-clockwise around its interior.
Then one can find four vertices of a square on $\gamma$ and order them such that they lie counter-clockwise on $\gamma$ and also label the square counter-clockwise.
\end{corollary}
\begin{proof}
This can be achieved by restricting $Q$ in the above proof to the set $Q'$ of $4$-tuples $[x_1,x_2,x_3,x_4]$ in $P_4$ that label the vertices of a square $(\gamma(x_1),\ldots,\gamma(x_4))$ in counter-clockwise order.
Along a bordism the orientation of the parametrised square cannot change (here we take a bordism that is induced by a isotopy of the curve in the plane).
If $\gamma$ is an ellipse then it is clear that $Q'$ equals $Q$, thus it represents the generator in $\unorientedBordism_0(P_4^\circ/_{\ZZ/4})$.
\end{proof}

\subsection{New cases: A simple open class of curves} \label{secSquarePegAsimpleOpenCondition}

In this section we prove Theorem~\ref{thmSquarePegInAnnulus} from the introduction together with the following two versions, whose proofs are very similar. 
See Figures~\ref{thmSquarePegInSquare} and~\ref{thmSquarePegInTriangle} for examples.

\begin{theorem}
\label{thmSquarePegInSquare}
Let $A$ denote the area 
$[-3,3]^2\wo(-1,1)^2$.
Suppose that $\gamma:S^1\to A$ is a continuous closed curve in $A$ that represents the generator of $\pi_1(A)=\ZZ$.
Then $\gamma$ inscribes a square of side length at least $\sqrt{2}$.
\end{theorem}


\begin{theorem}
\label{thmSquarePegInTriangle}
Let $\Delta$ be an equilateral triangle in $\RR^2$ whose center point is the origin.
Let $A$ be the closure of $((1+\sqrt{3})\cdot\Delta)\wo\Delta$.
Suppose that $\gamma:S^1\to A$ is a continuous closed curve in $A$ that represents the generator of $\pi_1(A)=\ZZ$.
Then $\gamma$ inscribes a square of side length at least $2\sqrt{3}-3$.
\end{theorem}

\begin{figure}[htb]
  \centering
  \begin{minipage}[b]{0.45\textwidth}
	\centering
    \includegraphics[scale=0.75]{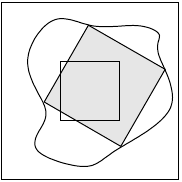}
    \caption{Example for Theorem~\ref{thmSquarePegInSquare}.}
    \label{figSquarePegInSquare}
  \end{minipage}
\quad
  \begin{minipage}[b]{0.45\textwidth}
    \centering
    \includegraphics[scale=0.75]{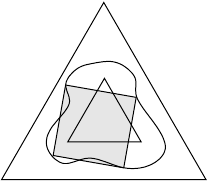}
	\caption{Example for Theorem~\ref{thmSquarePegInTriangle}.}
    \label{figSquarePegInTriangle}
  \end{minipage}
\end{figure}


It seems to be desirable to extend this method for much more general shapes in order to possibly prove the Square Peg Problem for all curves.

The proofs of Theorems~\ref{thmSquarePegInAnnulus},~\ref{thmSquarePegInSquare}, and~\ref{thmSquarePegInTriangle} follow from the following lemma.

\begin{lemma}
\label{lemSquareBordismLemma}
Let $A$ be an open subset of $\RR^2$.
Let $S_A$ be the set of $4$-tuples $(P_1,\ldots,P_4)\in A^4$ that form the vertices of a possibly degenerate square in counter-clockwise order.
Let $C$ be a connected component of an $\eps$-neighborhood of $S_A$ that does not contain degenerate squares, that is, points of the form $(P,P,P,P)$.
Let $\wt\gamma:S^1\to A$ be a generic curve that contains an odd number of squares in $C$.
Then every continuous curve $\gamma:S^1\to A$ that is homotopic to $\wt\gamma$ in $A$ contains a square in $C$ as well.
\end{lemma}

Here, by a generic curve $\wt\gamma$ we mean a curve such that the corresponding test map that measures squares in $C$ hits the test-space transversally.


\begin{proof}[Proof of Lemma~\ref{lemSquareBordismLemma}]
The proof is a simple bordism argument, similarly to Schnirelmann's proof above.
First assume that $\gamma$ is generic in the sense that its test map is transversal to the test-space.
When deforming $\gamma$ to $\wt\gamma$ within $A$, their inscribed squares change by a bordism, except that problems may appear close to degenerate squares.
However when we restrict the bordism to squares that are in the component~$C$, this problem disappears.
Thus, $\gamma$ and $\wt\gamma$ have modulo $2$ the same number of inscribed squares within~$C$, and by assumption this number is odd.

If $\gamma$ was not generic, there is several ways to make it generic.
Possibly the simplest one is to replace $\gamma$ by a slightly deformed copy within $A$ that is generic (this can be achieved by allowing deformations of $\gamma$ by local bumps and using the transversality theorem, as in Section~\ref{secTechnicalities}).
The previous argument then yields a square in $C$ that is almost inscribed in~$\gamma$.
Taking better and better approximations to $\gamma$ in the $C^0$-norm, and using the compactness of~$\gamma$, we obtain a square inscribed in~$\gamma$, which is non-degenerate as $C$ is bounded away from degenerate squares.
\end{proof}

\begin{proof}[Proof of Theorem~\ref{thmSquarePegInAnnulus}]
For $r\in [1,\infty]$, let $A_r=\{x\in\RR^2\st 1\leq ||x|| \leq r\}$, and consider its set of inscribed squares $S_{A_r}$ as defined in Lemma~\ref{lemSquareBordismLemma}.
For small enough~$r$, $S_{A_r}$ has two connected components, which we call the components of big squares and of small squares. 
The component of big squares deformation retracts $\ZZ/4$-equivariantly to the set of inscribed squares on the unit circle, and the component of small squares deformation retracts $\ZZ/4$-equivariantly to the set of degenerate squares in~$A_r$.
For large enough~$r$, $S_{A_r}$ becomes obviously connected.
To understand the phase transition, consider the path of squares $Q_\alpha$ ($\alpha\in [0, \frac{3}{4}\pi]$) with vertex set $\{e^{\pm i\alpha},e^{\pm i\alpha}+2\sin\alpha\}$. 
It connects the degenerate square $Q_0$ on the unit circle $S^1$ to a square $Q_{3\pi/4}$ that is inscribed in~$S^1$, such that for each~$\alpha$, two adjacent vertices of $Q_\alpha$ lie on~$S^1$ and the other two lie in the exterior of~$S^1$.
This is a path in $S_{A_\infty}$ between $Q_0$ and $Q_{3\pi/4}$ with the minimal possible $\max_\alpha \max_{P\in\textnormal{vert}(Q_\alpha)} ||P||$.
One computes that for the path $Q_\alpha$ this maximum vertex norm is $r_0:=1+\sqrt{2}$ and it is attained at $\alpha_0=\arctan(1+\sqrt{2})$.
This implies that for any $r<r_0$, $S_{A_r}$ has two connected components (as above), whereas for $r\geq r_0$, $S_{A_r}$ is connected. 
Moreover, $Q_{\alpha_0}$ has side length $(\sqrt{2}+2)^{1/2} > \sqrt{2}$, and we see that the smallest squares in the component of big squares of $S_{A,r}$ for $r<r_0$ are the ones inscribed in $S^1$, which have side length~$\sqrt{2}$.

To prove the theorem, we may assume that $\gamma$ is actually a curve in the interior of $A=A_{r_0}$.
The other cases follow by a limit argument: Approximate $\gamma$ by a sequence of curves in the interior, show that each approximating curve inscribes a square with side length at least $\sqrt{2}$, and take a convergent subsequence of such squares.
Thus we may replace $A$ by $A_{r^*}$ for some $r^*<r_0$, 
such that still $\gamma\subset A_{r^*}$.

Let $\wt\gamma$ be an ellipse in $A_{r^*}$ that represents a generator of $\pi_1(A_{r^*})$.
We apply Lemma~\ref{lemSquareBordismLemma} to $A_{r^*}$, $\gamma$ and~$\wt\gamma$, noting that the unique (generic) square inscribed in $\wt\gamma$ lies in the component of $S_{A_{r^*}}$ of big squares.
%
%
%
This proves the theorem.
\end{proof}

\begin{remark}
For any $r>1+\sqrt{2}$ it is open whether Theorem~\ref{thmSquarePegInAnnulus} holds for the annulus~$A_r$.
\end{remark}

The proofs of Theorems~\ref{thmSquarePegInSquare} and~\ref{thmSquarePegInTriangle} are analogous, using the fact that for a slightly smaller set $A$ the set of inscribed squares becomes disconnected.

\subsection{New cases: An open and dense class of curves} \label{secSquarePegMyProof}

First we will establish the main theorem of this section, which gives a quite general but somewhat technical condition for the existence of inscribed squares.
Then we deduce four less and less technical corollaries that demonstrate its scope. 

Let $\gamma:S^1\to\RR^2$ be a simple closed curve (that is, injective and continuous). We need some preparation. Let $f:P_4\to\RR^6$
be the corresponding test map that measure the four edges and two diagonals, which was defined in equation \eqref{eqTestMapSquarePeg} in Section \ref{secShnirelmansProof}.
For a path $\omega:S^1\to (S^1)^2\wo\diag_{(S^1)^2}$, $\omega(t)=(\omega_1(t),\omega_4(t))$, we define
\[
P_4(\omega) := \{[\omega_1(t),x_2,x_3,\omega_4(t)]\in P_4^\circ\st t\in S^1,\ x_2,x_3\in S^1 \textnormal{ such that } (\omega_1(t),x_2,x_3,\omega_4(t))\in\varphi(P_4)\},
\]
which is the set of quadrilaterals that are counter-clockwise inscribed in $S^1$ where the first and last vertex are of the form $\omega_1(t)$ and $\omega_4(t)$ for $t\in S^1$.

$P_4(\omega)$ can be parametrised by $g:S^1\times\simplex^2\to P_4(\omega)$, where $S^1$ parametrises $\omega$ and $\simplex^2$ the three arc lengths between the points $\omega_1(t)$, $x_2$, $x_3$ and~$\omega_4(t)$.
The map $g$ is injective if and only if $\omega$ is. 

\begin{definition}
We call an inscribed quadrilateral in $\gamma$ given by $[x_1,x_2,x_3,x_4]$ a \textdef{special trapezoid} if 
\begin{equation}
\label{eqDefSpecialTrapezoid}
f([x_1,x_2,x_3,x_4])=(a,a,a,b,e,e)\textnormal{ with }a>b\textnormal{, for some reals }a,b,e.
\end{equation}
The \textdef{size} of a special quadrilateral $[x_1,x_2,x_3,x_4]$ is the normalised arc length $x_4-x_1$.
\end{definition}

\begin{figure}[h]
\centering 
\input{Figures/SquarePegSmallSpecialTrapezoid3.pspdftex}
\caption{Example of a special trapezoid of size $\eps$.}
\label{figSpecialTrapezoid}
\end{figure}

Let $S\subset P_4$ denote the set of all special trapezoids.
Define $i(S,P_4(\omega))\in\FF_2$ as the mod-$2$ intersection number of $f\circ g$ and 
\begin{equation}
\label{eqV}
V:=\{(a,a,a,b,e,e)\in\RR^6\st a\geq b\}
\end{equation}
in $\RR^6$.
This is only well-defined if $f(g(S^1\times \bd \simplex^2))\cap V=\emptyset$ and $\im (f\circ g)\cap \bd V=\emptyset$.
The first requirement is trivially fulfilled because a quadrilateral $[x_1,x_2,x_3,x_4]$ with $x_i=x_{i+1}$ for some $i\in\{1,2,3\}$ cannot be special as $a>0$.
The second requirement holds if and only if $P_4(\omega)$ parametrised no inscribed quadrilateral $\gamma$ that is a square.
The map $f\circ g$ can now be deformed by a homotopy relative to $S^1\times\bd\simplex^2$, such that at no time it intersects the boundary of~$V$, and such that it becomes transversal to~$V$.
The intersection number then counts the preimages of $V$ under $f\circ g$ modulo~$2$. 

\begin{theorem}
\label{thmMySquarePeg}
Suppose there is a path $\omega: S^1\to (S^1)^2\wo\diag_{(S^1)^2}\homeo P_2^\circ$, $\omega(t)=(\omega_1(t),\omega_4(t))$, that represents a generator in $\pi_1((S^1)^2\wo\diag_{(S^1)^2})\iso\pi_1(S^1)\iso\ZZ$.
If $\gamma$ does not inscribe a square then the mod-$2$ intersection number $i(S,P_4(\omega))$ is well-defined and it equals~$1$.
\end{theorem}

The mod-$2$ intersection number will be described in the proof, see Section~\ref{secMySquarePegProof}.
First let us state and prove a few less technical corollaries.

\begin{corollary}
\label{cor1MySquarePeg}
Suppose there is a path $\omega: S^1\to (S^1)^2\wo\diag_{(S^1)^2}=P_2^\circ$, $\omega(t)=(\omega_1(t),\omega_4(t))$, that represents a generator in $\pi_1((S^1)^2\wo\diag_{(S^1)^2})\iso\pi_1(S^1)\iso\ZZ$.
If $P_4(\omega)\cap S=\emptyset$, then $\gamma$ has an inscribed square.
\end{corollary}
\begin{proof}
The mod-2 intersection number in Theorem \ref{thmMySquarePeg} is here trivially zero.
\end{proof}

\begin{corollary}
\label{cor2MySquarePeg}
Suppose there is an $\eps\in (0,1)$, such that $\gamma$ inscribes no (or generically an even number of) special trapezoid of size $\eps$.
Then $\gamma$ has an inscribed square.
\end{corollary}

\begin{proof}
Use Theorem \ref{thmMySquarePeg} with $\omega(t):=(t,t+\eps)$.
\end{proof}

\begin{remark}
The following piece of intuition might be useful.
Let $0<\eps<1/2$ and let us call an inscribed square $[x_1,x_2,x_3,x_4]$ of size larger than $\eps$ if $\{x_1,x_2,x_3,x_4\} \subset S^1$ does not lie in an arc of normalised length~$\eps$. 
Suppose everything is generic, and we change $\eps$ continuously or $\gamma$ by an ambient isotopy.
Then the set of inscribed special trapezoids of size $\eps$ changes by an unoriented bordism, except at points when $a$ becomes equal to $b$ (compare with Figure~\ref{figSpecialTrapezoid}), at which point an inscribed square of size larger than $\eps$ (with vertices cyclically on $\gamma$) will appear or disappear and the parity of the number of special trapezoids of size $\eps$ changes as well. 

We remark that Corollary~\ref{cor2MySquarePeg} holds indeed for any $\eps$, not only small ones.
\end{remark}

\begin{corollary}
\label{corSquarePegMainTheorem}
The class of curves in Corollary~\ref{cor2MySquarePeg} is an open and dense neighbourhood of Stromquist's locally monotone curves in the space of all injective maps $S^1\to\RR^2$ with respect to the $C^0$-topology.
In this sense, any generic Jordan curve has an inscribed square.
\end{corollary}

This implies Theorem~\ref{thmSquarePegMainTheorem} from the introduction.

\begin{proof}
The class of $\gamma$ without inscribed special trapezoids of size $\eps$ is open, even for fixed $\eps$.
By compactness, for any locally monotone curve $\gamma$, there exists an $\eps$ such that for any $x\in S^1$, there exists a linear function $\alpha:\RR^2\to\RR$ such that $\alpha\circ \gamma|_{[x,x+\eps]}$ is strictly increasing.
For this~$\eps$, $\gamma$ inscribes no special trapezoid of size~$\eps$, because a special trapezoid $ABCD$ with short edge $AB$ cannot be orthogonally projected to a line such that $A$ and $B$ project to the endpoints of the projection's image.
The density statement is obvious, since already the class of smooth curves is dense.
\end{proof}

\begin{corollary}
\label{cor4MySquarePeg}
Let $\gamma_0:S^1\incl\RR^2$ be a $C^2$-embedding with curvature bounded above by~$k$.
Then any continuously embedded $\gamma: S^1\to\RR^2$ in the $\frac{1}{4k}$-neighbourhood of $\gamma_0$ with respect to the $C^0$-metric has an inscribed square.
\end{corollary}

\begin{proof}
We may and do assume that $\gamma_0$ is parametrised by arc-length and that $k=1$. Let $\delta=\frac{1}{4}$.
The following calculation shows that $\gamma$ does not inscribe a special trapezoid of size $\eps=2\arccos(\frac{1}{3})$.

Assume $ABCD$ was such a hypothetical special trapezoid with parameters $x_1,x_2,x_3,x_4$ that satisfy~\eqref{eqDefSpecialTrapezoid} and $x_4-x_1=\eps$. They also parametrise an inscribed quadrilateral $A_0B_0C_0D_0$ in $\gamma_0$. 
We have $a\geq b=||A-D|| \geq ||A_0-D_0||-2\delta \geq 2\sin\frac{\eps}{2} - 2\delta =: d_1$, where the first inequality is the triangle inequality, and the latter one uses the given curvature bound and the formula for the length of a chord in the unit circle with central angle~$\eps$.
On the other hand, for some $1\leq i \leq 3$, $x_{i+1}-x_i\leq \frac{\eps}{3}$, say for $i=1$ (the other cases are analogous).
Then $a = ||B-A|| \leq ||B_0-A_0||+2\delta \leq \frac{\eps}{3}+2\delta := d_2$.
Together, both bounds yield $d_2\leq d_1$, which is impossible for the chosen values of $\delta$ and~$\eps$.


The result follows from Corollary~\ref{cor2MySquarePeg}.
\end{proof}

\begin{remark}
We remark that in Corollary~\ref{cor4MySquarePeg}, $\delta=\frac{1}{4k}$ is not best possible. It could be replaced by $\frac{1}{3k}$, as a more detailed analysis shows that there will be no inscribed special trapezoids of size $\eps=\frac{\pi}{2k}$. 
It is to check that in any arc $a$ of length $\eps$ there is no special trapezoid $ABCD$ with $AD$ as the shortest edge such that $A$ is $\frac{1}{3k}$-close to one end point of $a$, $D$ is $\frac{1}{3k}$-close to the other end point, and $B$ and $C$ are $\frac{1}{3k}$-close to $a$, see Figure~\ref{figSausage}.
\begin{figure}[h]
\centering 
\input{Figures/CurvatureBound.pspdftex}
\caption{$\frac{1}{3k}$-neighborhood of $a$. In this particular figure, $a$ has constant curvature $k$.}
\label{figSausage}
\end{figure}
\noindent
On the other hand, for example $\delta=\frac{1}{2k}$ is out of reach with our method.
\end{remark}

\subsubsection{Proof of Theorem~\ref{thmMySquarePeg}}
\label{secMySquarePegProof}
The proof is based on equivariant obstruction theory. This was first used in connection to the Square Peg Problem by \vrecica\ and \zivaljevic~\cite{VrZi08fultonMacPhersonCompCyclohedraAndPolygonalPegProblem}, who applied it to the Fulton--MacPherson compactification of the configuration space~$P_4^\circ$.
Here, we will apply it to suitable closed subsets of~$P_4^\circ$. 

We will need the following instance of basic relative obstruction theory.
Consider a rank $k$ vector bundle $\varphi:E\to B$ over a connected $k$-dimensional manifold $B$ with non-empty connected boundary~$\bd B$.
Let $0_\varphi\subset E$ denote the zero section.
Given a nowhere vanishing section $u:\bd B\to E$ over the boundary, there is a relative Stiefel--Whitney class $w_k(u)\in H^k(B,\bd B;\FF_2)\iso H_0(B;\FF_2) = \FF_2$, which is the primary obstruction for the existence of a nowhere vanishing extension of $u$ to all of~$B$ modulo $2$.
If $s$ is an arbitrary extension of $u$ to all of $B$, then the parity of the generic number of zeros of $s$ coincides with the Poincar\'e dual of $w_k(u)$ in~$\FF_2$.
If $u_1,u_2$ are two nowhere vanishing sections of $\varphi|_{\bd B}$ then the primary obstruction to the existence of a fiberwise homotopy between $u_1$ and $u_2$ is the so-called difference cocycle $d(u_1,u_2)\in H^{k-1}(\bd B;\FF_2)\iso H_0(\bd B;\FF_2)=\FF_2$.
Moreover, if $\delta$ denotes the connecting homomorphism $H^{k-1}(\bd B;\FF_2)\to H^k(B,\bd B;\FF_2)$, which is an isomorphism in our setting since $H^k(B;\FF_2)=0$, then $\delta(d(u_1,u_2))=w_k(u_1)-w_k(u_2)$.
Moreover, choose a third nowhere vanishing section $t$ over~$\bd B$.
Considering $t,u_i$ as sections in the associated sphere bundle of $\varphi_{\bd B}$, we can define mod-2 intersection numbers $i(t,u_i)$, $i=1,2$.
Then the Poincar\'e dual of $d(u_1,u_2)$ equals $i(t,u_1)-i(t,u_2)$.
Therefore, considered as elements of $\FF_2$, $w_k(u_1)-w_k(u_2) = i(t,u_1)-i(t,u_2)$.

\begin{proof}[Proof of Theorem~\ref{thmMySquarePeg}]
First we prove the theorem in the special case when $\omega$ is the path
$\omega_\eps(t):=(t,t+\eps)$ for some fixed $0<\eps<1/2$; we may choose $\eps=1/4$.
Then $P_4(\omega)\cap S$ is the set of inscribed special trapezoids of size $\eps$.

Let 
\[
P_4^\eps:=\{(x;t_1,\ldots,t_4)\in P_4\st t_1,\ldots,t_4\leq 1-\eps\} = S^1\times T^\eps,
\]
where $T^\eps:=\{(t_1,\ldots,t_4)\in \simplex^3\st t_1,\ldots,t_4\leq 1-\eps\}$ is a truncated tetrahedron.
Let $T_1',\ldots,T_4'$ be the triangular facets of $T^\eps$, and $H_1',\ldots,H_4'$ the hexagonal facets.
That is, $T_i'=\{t\in T^\eps\st t_i=1-\eps\}$ and $H_i'=\{t\in T^\eps\st t_i=0\}$.
Moreover, define $T_i := S^1\times T_i'$ and $H_i := S^1\times H_i'$ for $1\leq i \leq 4$.
Note that for our choice of $\omega$, $P_4(\omega)=T_4$.
Moreover, $\bd P_4^\eps = \bigcup_i T_i\ \cup\ \bigcup H_i$.

Now assume that $\bd P_4^\eps$ parametrises no inscribed square (otherwise we are done).
Then the parity of the generic number of inscribed squares up to $\ZZ/4$-symmetry that are parametrised by $B:=P_4^\eps/_{\ZZ/4}$ is well-defined, call it~$n^\eps(\gamma)$.

Let $\varphi:E\to B$ denote the canonical projection from $E:=P_4^\eps\times_{\ZZ/4} (\RR^6/\Delta)$ to~$B$.
This is a rank~$4$ vector bundle with fiber $\RR^6/\Delta$.
Let $0_\varphi$ denote its zero-section.
The test map $f$ induces a section $s_1:B\to E$ and $s_1^{-1}(0_\varphi)$ is the set of all inscribed squares up to $\ZZ/4$-symmetry that are parametrised by~$B$, and put $u_1:=s_1|_{\bd B}$.
Moreover $n^\eps(\gamma)$, the generic parity of $s_1^{-1}(0_\varphi)$ (assuming that $u_1$ is nowhere zero), equals the Poincar\'e dual of the top relative Stiefel--Whitney class $w_4(u_1)\in H^4(B,\bd B;\FF_2)\iso \FF_2$.

Similarly let $s_2:B \to E$ be the analog section that is induced from an arc-length parametrised ellipse $\gamma_0$ instead of $\gamma$, and let $u_2:=s_2|_{\bd B}$.

By the above discussion, $n^\eps(\gamma)-n^\eps(\gamma_0) = i(t,u_1)-i(t,u_2)$ (as elements in $\FF_2$) for any nowhere vanishing section $t$ of $\varphi|_{\bd B}$.
A nowhere vanishing section $t$ of $\varphi|_{\bd B}$ is the same thing as a $\ZZ/4$-equivariant map $t':\bd P_4^\eps \to (\RR^6/\Delta)\wo 0$.
We construct a suitable $t$ via $t'$ as follows.
For any $p=(x;t_1,\ldots,t_4)\in T_4\subset \bd P_4^\eps$ define $t'(p)$ as $(1,1,1,0,0,0)+\Delta$.
Let $R\subset\RR^6/\Delta$ be the ray from $0$ through this vector, that is, $R/\Delta = V/\Delta$, where $V$ was defined in~\eqref{eqV}.
Then $f(p)\in R$ if and only if $p$ parametrises a special trapezoid or a square.
Extend $t'$ $\ZZ/4$-equivariantly to $\bigcup_i T_i$. 

Now extend $t'$ furthermore continuously on $H_4$ such that for all $p=(x;t_1,\ldots,t_4)\in H_4$, $t'(p)=(a_1,a_2,a_3,a_4,0,0)+\Delta$ with $a_4>\min(a_1,a_2,a_3)$, and such that $t'$ is equivariantly extendable to $\bigcup_i H_i$. 
An explicit way to do this is to define
\[
\textstyle
t'(p) = \Big(r(t_1),r(t_2),r(t_3),1+\prod\limits_{i=1}^3 t_i(1-\eps-t_i),0,0\Big)+\Delta, \quad \textnormal {for }p=(x;t_1,\ldots,t_4)\in H_4,
\]
where $r:[0,1-\eps]\to [0,1]$ is defined by $r(x)=1$ for $x
\in[0,1/2]$, $r(1-\eps)=0$, and linearly between $1/2$ and $1-\eps$.

Any inscribed quadrilateral that is parametrised by $p\in H_4$ has a degenerate third edge.
Therefore $\RR_{\geq 0}\cdot f(p)$ will never coincide with $\RR_{\geq 0}\cdot t'(p)$ for $p\in H_4$.
Equivalently, $\RR_{\geq 0}\cdot s_1([p])$ will never coincide with $\RR_{\geq 0}\cdot t([p])$ for $p\in H_4$, where $[p]$ denotes the image of $p$ in $B=P_4^\eps/_{\ZZ/4}$.

Thus $i(t,u_1)$ is precisely the parity of the generic number of inscribed special trapezoids of size $\eps$ for~$\gamma$, and $i(t,u_2)$ is the same for the ellipse.

By inspecting the ellipse, which is particularly easy for small~$\eps$, we get $n^\eps(\gamma_0)=1$ and $i(t,u_2)=0$, hence $n^\eps(\gamma_0)-i(t,u_1)=1$ in~$\FF_2$.
Therefore for $\gamma$ we get $n^\eps(\gamma) = 1 + i(t,u_1)$ in~$\FF_2$, which proves the theorem for the above~$\omega_\eps$.

\medskip

For general $\omega$, we reduce to the previous argument for $\omega_\eps$ as follows.
By possibly reversing the orientation of $\omega$, we may assume that $\omega$ and $\omega_\eps$ are homotopic within $(S^1)^2\wo\Delta_{(S^1)^2}$ via some homotopy~$\Omega_\tau$, $\tau\in [0,1]$, with $\Omega_0 = \omega_\eps$ and $\Omega_1 = \omega$.
Using $\Omega_\tau$ we will construct a $\ZZ/4$-map $h:P_4^\eps \to P_4^\circ$ such that for each $1\leq i\leq 4$, $H_i$ is sent into $\{(x;t_1,\ldots,t_4)\st t_i=0\}$, and such that $T_4\homeo S^1\times\simplex^2$ is sent homeomorphically to~$P_4(\omega)$. 
Explicitly, we may define $h$ as follows. 
Choose $\eps' \in (\eps,1/2)$, e.g.~$\eps'=1/3$ if $\eps=1/4$.
On the subset~$P_4^{\eps'}$, we let $h$ be the identity.
The complement $P_4^{\eps}\wo P_4^{\eps'}$ has four connected components, one of which is the set of all $p = (x;t_1,\ldots,t_4) \in P_4^\eps$ that satisfy $t_4 > 1-\eps'$.
At such $p$, we define 
$\tau^{(p)} := ((1-\eps)-t_4)/(\eps'-\eps)$, 
$\omega^{(p)} = (\omega^{(p)}_1,\omega^{(p)}_4):= \Omega_{\tau^{(p)}}$, 
$c^{(p)} := \tau^{(p)} + (1-\tau^{(p)})(\omega^{(p)}_4-\omega^{(p)}_1)/\eps$, and
\[
h(p) = (\omega_1^{(p)}(x);c^{(p)}t_1,c^{(p)}t_2,c^{(p)}t_3,1-c^{(p)} (1-t_4)).
\]
We extend $h$ to the other three connected components of $P_4^{\eps}\wo P_4^{\eps'}$ by $\ZZ/4$-equivariance. This makes $h$ well-defined, continuous, and it satisfies the mentioned properties.
Let $s_1': B \to E$ be the section induced by $f\circ h$, and let $u_1' = s_1'|_{\bd B}$.
As with $\omega_\eps$, we apply an analogous obstruction theory argument, this time to $s_1'$ instead of~$s_1$, to obtain $w_4(u_1') = 1 + i(t,u_1') \in \FF_2$.
The difference is that now, (i) $i(t,u_1')$ counts the parity of the generic number of inscribed special trapezoids in $P_4(\omega)$, instead of those of size~$\eps$, and (ii) the \poincare{} dual of the relative Stiefel--Whitney class $w_4(u_1')$ counts via $s_1'$ the generic number of $\ZZ/4$-orbits $[p]\in B$ such that $h(p)$ parametrises an inscribed square in~$\gamma$. 
This proves the theorem.
\end{proof}

\subsubsection{Remarks} \label{secMySquarePegRemarks}

\begin{enumerate}[1.]
\item The proof of Theorem~\ref{thmMySquarePeg} gives slightly more:
Let $q_\eps$ be the parity of the generic number of inscribed squares up to symmetry whose parametrisation $(x;t_0,\ldots,t_3)\in P_4$ satisfies $\sharp\{i\st t_i\leq 1-\eps\}=0\mod 2$.
And let $s_\eps$ be the parity of the generic number of inscribed special trapezoids of size~$\eps$.
If one of $q_\eps$, $s_\eps$ is well-defined then both of them are, and $q_\eps\neq s_\eps\mod 2$.
An analog statement holds for arbitrary~$\omega$.

\item 
%
$P_4(\omega)$ can be regarded as a ``membrane'', which separates $P_4$ into two components if $\omega$ is injective.
If $\gamma$ circumscribes no square then there is an \emph{odd} number of paths in $S$ that pass through $P_4(\omega)$ and approach the boundary at $S^1\times e_3$, $e_3$ being the one vertex of~$\simplex^3$.
These paths might look very chaotic close to the boundary.
On the other side of the membrane $P_4(\omega)$, this odd number of paths cannot all end in each other.
One of them has to end somewhere else.
It might end suddenly in $P_4^\circ$, which means that it found a square, or it might end somewhere else at~$\bd P_4^\circ$.
The latter is (unfortunately) possible:

\begin{figure}[h]
\centering 
\input{Figures/Spirale3.pspdftex}
\end{figure}

The drawn path of special trapezoids starts in the middle of the spiral at $S^1\times e_3$ with a quadrilateral that is degenerate to a point, and it stops when $x_1$ and $x_4$ meet again, $x_4-x_1=1$. 

\item The original proof of Theorem~\ref{thmMySquarePeg} from~
\cite{Mat11phd}
 is slightly more elementary (but admittedly more difficult to digest).

\item The corollaries are sometimes good for proving the existence of a square, if the curve is piecewise $C^1$ but has cusps (points in which the tangent vector changes the direction). This however works not in a large generality as the previous example shows.

\item Theorem~\ref{thmMySquarePeg} and its proof deal with the curve \emph{intrinsically}, since the only datum of $\gamma$ we used is the distances between points on $\gamma$.
If we define a square in a metric space $(X,d)$ to be a 4-tuple $(x_0,\ldots,x_3)\in X^4$ such that $d(x_0,x_1)=d(x_1,x_2)=d(x_2,x_3)=d(x_3,x_0)$ and $d(x_0,x_2)=d(x_1,x_3)$, then the theorem also holds for curves $\gamma:S^1\to X$.
More generally, $X$ does not need to fulfill the triangle inequality.
In other words, we do not need an embedded curve but a distance defining function $d:S^1\times S^1\to\RR$ that is continuous, positive definite, and symmetric.
\end{enumerate}

\section{Immersed curves} \label{secImmersedCurves}

\subsection{Squares inscribed in immersed curves}

Toeplitz's conjecture is about inscribed squares on simple closed curves in the plane. 
There are plenty of ways to generalize this problem.
In this section we study what we can say about the number of inscribed squares if we omit the requirement that $\gamma$ has to be injective.
In this setting it makes sense to more generally allow $\gamma$ to be a finite union of curves, $\gamma: X\to\RR^2$, $X=\coprod_{i=1}^n S^1$.
There are several kinds of degenerate squares, which we have to deal with in that case. 
To be able to count inscribed squares in a stable manner, we will only consider generic smooth curves and we will not count degenerate squares.
Here, \emph{generic} means that:
(i) self-intersections of $\gamma$ are transversal with a non-orthogonal intersection angle, 
(ii) the test map~\eqref{eqTestMapForImmersions} is transversal to the test-space~$V$, and 
(iii) no non-degenerate inscribed square has a vertex at a self-intersection point of~$\gamma$.
If any of these condition fails, the number of inscribed squares might not be stable under small ambient isotopies of~$\gamma$.

Given $\gamma: X\to\RR^2$, we construct a test map 
\begin{equation}
\label{eqTestMapForImmersions}
t:X^4\wo D \to \RR^6\supset \Delta
\end{equation}
with the same pointwise formula as in~\eqref{eqTestMapSquarePeg}, $\Delta = \Delta_{\RR^4}\times \Delta_{\RR^2}$ as before, and where $D = \{x\in X^4\st \gamma(x_1)=\gamma(x_2)=\gamma(x_3)=\gamma(x_4)\} = (\gamma^4)^{-1}(\Delta_{(\RR^2)^4})$ parametrises the set of degenerate inscribed squares.



We call the self-intersections of $\gamma$ crossings.
These crossings together with the connecting arcs of $\gamma$ form a planar graph, which is $4$-regular.
Hence its dual graph is bipartite, which means that we can (uniquely) color the components of the complement of $\gamma$ in black and white such that adjacent components obtain different colors and such that the unbounded component is white.
This is sometimes called a chequerboard coloring of the complement of $\gamma$,
%
see~Figure~\ref{figChequerBoardColoring}. 
%
%
%
Let~$b(\gamma)$ be the number of black components.
We say that a crossing is \textdef{fat} if the black angles at this crossing are larger than~$90^\circ$. 
The fat crossings in Figure~\ref{figChequerBoardColoring} are marked with black dots.
Let~$f(\gamma)$ be the number of fat crossings.

\begin{figure}[tbh]
\centering
\input{Figures/ChequerBoardColoring.pspdftex}
\caption{Chequerboard coloring associated to $\gamma$. Dots mark the fat crossings.}
\label{figChequerBoardColoring}
\end{figure}

\begin{theorem}
\label{thmSquaresOnImmersedCurves}
Suppose that~$\gamma$ is an immersion of a closed $1$-manifold in the plane that is generic in the above sense. 
Then the number of non-degenerate squares inscribed in $\gamma$ is congruent to $b(\gamma)+f(\gamma)$ modulo~$2$.
\end{theorem}

\begin{proof}
By genericity of the curve, no inscribed square will have a vertex at a crossing.
Thus at each crossing $c$ we can find a neighborhood $N_c \subset \RR^2$ of $c$ that contains no vertex of a non-degenerate inscribed square of $\gamma$ and in which $\gamma$ is $C^2$-close to two intersecting straight line segments.
We may deform $\gamma$ slightly (in the $C^1$-sense) within each $N_c$, such that in a suitably smaller neighborhood $N_c'$ of $c$, $\gamma$ becomes the intersection of two such straight line segments with the same intersection angle.
If $N_c'$ is small enough, the set of non-degenerate inscribed squares stays invariant: Inscribed squares with parameter $4$-tuple close to $D$ cannot appear due to the smoothness of~$\gamma$; others don't appear due to the continuity of~$t$ and the smallness of~$N_c'$.

%
Next, at each crossing $c$ we perform a surgery as in Figure~\ref{figSmootheningACrossing} in a sufficiently small neighborhood $N''_c\subset N'_c$ of $c$, smoothening the crossing with two convex arcs in such a way that the white angles open up.
%
In this way, all white components merge into one unbounded component, and the number of black component stays invariant.

We need to study how the surgery changes the number of inscribed squares:
Non-degenerate squares that were inscribed before the surgery stay inscribed.
New inscribed squares must have at least one vertex in $N''_c$ for some $c$.
We argue that in fact all four vertices lie in $N'_c$: If not, then the $4$-tuple parametrising this square is bounded away from $D$, and hence by continuity of $t$ we may assume that all $N''_c$ were chosen small enough so that such a new square could not have appeared after the surgery.
It remains to consider newly inscribed squares with all four vertices in some $N'_c$.
Such a square appears in $N'_c$ if and only if the black angle at $c$ is larger than $90^\circ$, in which case the square is unique (see Figure~\ref{figSmootheningACrossing}): 
If the two arcs are arcs of an ellipse, this is an elementary but technical calculation, which we omit.

Thus, during the surgery, the number of inscribed squares increases by $f(\gamma)$. 
The new curve consists of $b(\gamma)$ separated simple closed curves. 
We can deform them by an ambient isotopy such that they become $b(\gamma)$ sufficiently small ellipses with mid-points aligned on a fixed line with sufficiently large distance between them, such that there is trivially no inscribed square with vertices in more than one component. 
Therefore the resulting union of ellipses inscribes exactly $b(\gamma)$ squares. 
Using a bordism argument as in Section~\ref{secShnirelmansProof}, the parity of the number of inscribes squares did not change during the isotopy. 
Since every ellipse inscribes exactly one square, this finishes the proof.
\end{proof}

\begin{figure}[tbh]
\centering
\input{Figures/SmootheningACrossing3.pspdftex}
\caption{When we smoothen a non-orthogonal straight line crossing then a new square appears if and only if we opened the smaller angle.}
\label{figSmootheningACrossing}
\end{figure}

\subsection{Rectangles inscribed in immersed curves}

The analog theorem for rectangles of prescribed aspect ratio $0<r<1$ that are inscribed in an immersed closed $1$-manifold is slightly different.
Let $0<\alpha(r)<\pi/2$ be the intersection angle of the two diagonals of a rectangle with aspect ratio~$r$.
Let $\gamma$ be a generic immersion of a finite union of circles in the plane, and consider again the chequerboard coloring from above.
Here, generic means that 
(i) all crossings are transversal with intersection angles are neither $\alpha(r)$ nor $\pi-\alpha(r)$, 
(ii) the test map $t_r:X^4\wo D \to \RR^4$ with the same formula as in~\eqref{eqTestMapRectangles} is transversal to~$(0,0,r)$ (i.e.\ all non-degenerate inscribed rectangles with aspect ratio $r$ are generic), 
and (iii) no (non-degenerate) inscribed rectangle with aspect ratio $r$ has a vertex at a self-intersection point of~$\gamma$.

We call a crossing of $\gamma$ \textdef{$\alpha$-orthogonal}, if its angle lies in the interval $(\alpha,\pi-\alpha)$.
Let $o(\gamma,r)$ denote the number of $\alpha(r)$-orthogonal crossings.

\begin{theorem}
\label{thmRectanglesOnImmersedCurves}
Let $0<r<1$.
Suppose that $\gamma$ is a generic immersion of a closed $1$-manifold in the plane. 
Then the number of non-degenerate rectangles with aspect ratio $r$ inscribed in $\gamma$ is congruent to $b(\gamma)+o(\gamma,r)$ modulo $2$.
\end{theorem}

\begin{proof}
The proof is analogous to the one of Theorem~\ref{thmSquaresOnImmersedCurves}. 
The only difference is that when we smoothen the crossing, then zero, one, or two new rectangles will appear, depending on whether the angle $\beta$ that we smoothen satisfies $\beta<\alpha$, $\alpha<\beta<\pi-\alpha$, or $\pi-\alpha<\beta$; compare with Figure~\ref{figSmootheningACrossingRectangle}.
\end{proof}

\begin{figure}[tbh]
\centering
\input{Figures/SmootheningACrossingRectangle.pspdftex}
\caption{Smoothening a crossing changes the parity of the number of inscribed rectangles with aspect ratio~$r$ if and only if the crossing is $\alpha$-orthogonal.}
\label{figSmootheningACrossingRectangle}
\end{figure}

\section{Inscribed rectangles} \label{secRectanglesOnCurves}

The Rectangular Peg Problem, Theorem~\ref{thmRectangularPegProblem}, was a very challenging and from the author's point of view the most beautiful open problem%
\footnote{Update: see Section~\ref{secRecentProgress}.} 
in this area of inscribing and circumscribing problems.



Griffiths \cite{Gri91topSquarePeg} gave a proof, however there are unfortunately some errors in his computation concerning orientations (see \cite[Chap. III.7]{Mat08diplomaThesis} for details). 
Hence the problem is open.

In Section~\ref{secRectanglesCSTMfails} we show that the standard topological approach, the \emph{configuration space/test map scheme}, fails to prove the Rectangular Peg Problem since the test map in question exists.
Then we prove Theorem~\ref{thmRectangularPegForSqrt3} under some technical assumptions concerning trans\-ver\-sa\-li\-ty; see Section~\ref{secRectanglesProofOfMainTheorem}.
We show in Section~\ref{secTechnicalities} that these assumptions can be made.
These technicalities seem not to be obvious in advance for two reasons:
The natural group action on one solution manifold (namely $P$) is in general not free; and transversality has to be achieved for several maps simultaneously since we need to relate solution manifolds of different maps in the proof.

\subsection{The test map exists}
\label{secRectanglesCSTMfails}

In this section we show that a standard topological approach, called configuration space/test map method, does not work to prove the Rectangular Peg Problem~\ref{thmRectangularPegForSqrt3}. 

Assume we are given a smooth simple closed planar curve $\gamma: S^1\incl\RR^2$. 
As above let $P_4^\circ\subset (S^1)^4$ be the set of four pairwise distinct points on the circle that lie counter-clockwise on it. 
Then $P_4^\circ\homeo S^1\times (\simplex^3)^\circ$, where $(\simplex^3)^\circ$ denotes the interiour of the $3$-simplex. 
We construct from $\gamma$ a natural test map,
\begin{equation}
\label{eqTestMapRectangles}
\begin{aligned}
t:&~P_4^\circ&\To &\ \RR^2\times\RR\times \RR,\\
&~[x_1,x_2,x_3,x_4]&\longmapsto &\ (v,\ell,a),
\end{aligned}
\end{equation}
where $v$ is the difference between the midpoints of the diagonals in the quadrilateral with vertices  $\gamma(x_1),\ldots,\gamma(x_4)$; $\ell$ is the difference of the length of these diagonals, and $a$ is the aspect ratio.

We let $\ZZ/2=\{0,\eps^2\}\subset\ZZ/4$ act on $P_4^\circ$ by $\eps^2\cdot[x_1,x_2,x_3,x_4]=[x_3,x_4,x_1,x_2]$. 
The map $t$ is then $\ZZ/2$-equivariant with respect to the corresponding group action on~$\RR^4$.
Since $\gamma$ is smooth, there is an $\eps>0$, such that $t$ maps no point of $B:=U_{[\eps]}(\bd(P_4^\circ))\cap P_4^\circ$ to zero, where $U_{[\eps]}$ denotes the closed $\eps$-neighborhood and $\bd P_4^\circ$ the topological boundary of $P_4^\circ\subset P_4$. 
The map $t|_B:B\to\RR^4\wo\{0\}$ is uniquely given up to $\ZZ/2$-homotopy. 
$R:=t^{-1}(0,0,r)$ is the set of rectangles of aspect ratio $r$ whose vertices lie counter-clockwise on~$\gamma$. 
$R$ is generically a zero-dimensional free $\ZZ/2$-manifold. 
Using the preimage orientation for $R$, $\ZZ/2$ acts orientation preserving on~$R$. 
Therefore $R$ determines an element $[R]$ in the oriented zero-dimensional bordism group $\Omega_0(P_4^\circ/_{\ZZ/2})\iso\ZZ$, which is 
the primary obstruction for extending $t|_U:U\to \RR^4\wo\{0\}$ to a map $P_4^\circ\to\RR^4\wo\{0\}$. 
If $\gamma$ is an ellipse then $R$ consists of two orbits, since an ellipse inscribes exactly two rectangles. 
A computation shows that their orientation is opposite, see~\cite{Mat09squarePegRepository}. %
Therefore $[R]=0$ and the obstruction class vanishes. 
Since this is the only obstruction, we can find a map $t':P_4^\circ\to\RR^4\wo\{0\}$ such that $t'|_B=t|_B$.
That is there is no purely topological argument that can show the existence of a rectangle of aspect ratio $r$ on $\gamma$, at least as long as we are not using more geometric information.

The smooth Square Peg Problem can be solved using this configuration space/test map scheme, since squares have more symmetry. 
Here the group of symmetry is $\ZZ/4$ and on an ellipse we find only one $\ZZ/4$-orbit of squares.

\subsection{Topological criteria} \label{secRectanglesTopologicalIntuition}


Above we saw that due to the lack enough symmetry, purely topological arguments will not work to prove the Rectangular Peg Problem.
But they give some intuition, here are two approaches. Assuming that Theorem~\ref{thmRectangularPegProblem} admitted a counter-example $(\gamma,r)$, both lemmas derive conclusions that seem to be unintuitive, but more geometric ideas are needed to yield a contradiction.

If $\alpha:S^1\to P_4^\circ$ is a one-parameter family of quadrilaterals, then we call $[p_1\circ\alpha]\in\pi_1(S^1)=\ZZ$ its winding number, where $p_1:P_4^\circ\to S^1$ denotes the projection to the first coordinate.

\begin{lemma}
\label{lemRectanglesOnCurvesUsingParallelograms}
Suppose there was a counter-example $(\gamma, r)$ for the (smooth) Rectangular Peg Problem. Then for all $\eps>0$, there is a $\ZZ/2$-invariant one-parameter family $S^1\to P_4^\circ$ of $\eps$-close parallelograms with aspect ratio in $[r-\eps,r+\eps]$ and with an odd winding number, such that during the whole one-parameter family one of the diagonals stays larger than the other one.
\end{lemma}

\begin{proof}



Given $\gamma$, $r$ and $\eps$ as in the lemma, let $\eps_n:=\eps/n$ for some $n\geq 1$.
Consider the $\ZZ/2$-equivariant test map
\[
\begin{aligned}
g:\ &P_4^\circ &\To&\ \ \RR^2\times\RR \\
&[x_1,x_2,x_3,x_4]&\longmapsto&\ \bigl((\gamma(x_1)+\gamma(x_3))-(\gamma(x_2)+\gamma(x_4)),\\
&&&\ \ (||\gamma(x_1)-\gamma(x_2)||+||\gamma(x_3)-\gamma(x_4)||)- \\
&&&\ \ r\cdot(||\gamma(x_2)-\gamma(x_3)||+||\gamma(x_4)-\gamma(x_1)||)\bigr).
\end{aligned}
\]
The preimage $P=(g|_{P_4^\circ})^{-1}(0)$ of the fixed point $0$ is the $\ZZ/2$-invariant set of inscribed parallelograms of aspect ratio~$r$ whose vertices lie cyclically on~$\gamma$.
As $\gamma$ is smooth, 
$\bd P_4^\circ$ has a closed neighborhood $N$ in $P_4$ that does not intersect~$P$.
As $\ZZ/2$ acts moreover freely on~$P_4^\circ$, we can deform $g$ by an equivariant $\eps_n$-homotopy relative to $N$ to a new map $g':P_4^\circ\to\RR^3$ that is transversal to~$0$.
Its preimage $P'=g'^{-1}(0)$ becomes a one-dimensional compact $\ZZ/2$-manifold of $\eps_n$-close inscribed parallelograms.
The bordism argument of Section~\ref{secShnirelmansProof} carries over to show that $[P'/_{\ZZ/2}]\in \unorientedBordism_1(P_4^\circ/_{\ZZ/2})=\unorientedBordism_1(S^1/_{\ZZ/2})=\ZZ/2$ represents the same element as $P$ does when $\gamma$ is the unit circle, and thus this element is the generator of~$\ZZ/2$.
Let $C$ be a component of $P'$.
If $C$ is not $\ZZ/2$-invariant, then $[(\ZZ/2\cdot C)/_{\ZZ/2}]$ is $0$ in $\unorientedBordism_1(P_4^\circ/_{\ZZ/2})$.
If $C$ is $\ZZ/2$-invariant, then $[C/_{\ZZ/2}]$ is the generator of $\unorientedBordism_1(P_4^\circ/_{\ZZ/2})$ if and only if the winding number of $C$ is odd.
Thus $P'$ has an odd number of $\ZZ/2$-invariant components with odd winding number.
Consider such a component~$C$.
If for some $n$, along the one-parameter family of $\eps_n$-close parallelograms parametrised by $C$, one diagonal stays always longer than the other, then the lemma is proved.

Otherwise for each $n$, by the intermediate value theorem at least one of these quadrilaterals must be an $\eps_n$-close rectangle~$R_n$.
As $\gamma$ is smooth, $(R_n)_n$ can be uniformly bounded away from $\bd P_4^\circ$.
Thus some convergent subsequence of $(R_n)_n$ converges to a non-degenerate inscribed rectangle of aspect ratio~$r$, which is a contradiction.
\end{proof}

\begin{remark}
In Lemma \ref{lemRectanglesOnCurvesUsingParallelograms}, instead of considering the set of parallelograms with aspect ratio~$r$, we might look as well on the set of parallelograms whose diagonals intersect in an angle $\alpha$, where $\alpha$ is the intersection angle of the diagonals in a rectangle of aspect ratio $r$. This gives an analogous lemma, which might be easier to deal with geometrically.
\end{remark}


\begin{lemma}
\label{lemRectanglesOnCurvesUsingRectangles}
Suppose there was a counter-example $(\gamma, r)$. Then for all $\eps>0$, there is a $\ZZ/4$-invariant one-parameter family $S^1\to P_4^\circ$ of $\eps$-close rectangles.
\end{lemma}
\begin{proof}
Let $f:P_4^\circ\toG{\ZZ/4} \RR^4\times \RR^2$ be the restricted map \eqref{eqTestMapSquarePeg} from Section \ref{secShnirelmansProof}, measuring the edges and diagonals of inscribed quadrilaterals. 

First, make $f$ $\ZZ/4$-equivariantly transversal to $\Delta=\diag_{\RR^4}\times\diag_{\RR^2}$ by a small $\delta$-homotopy, where $\delta$ is a positive function that decreases sufficiently fast near the boundary of $P_4^\circ$, and let $Q:=f^{-1}(\Delta)$ be the set of all squares (as measured with an error bounded by~$\delta$).
Then make $f$ $\ZZ/4$-equivariantly transversal to the $\ZZ/4$-invariant subspace $\Delta':=\{(a,b,a,b,e,e)\in\RR^4\times\RR^2\}$ by a small $\delta$-homotopy which leaves $Q$ fixed, and let $R:=f^{-1}(\Delta')$. 
If $\delta$ was chosen small enough in terms of a given~$\eps>0$, $R$ parametrises inscribed $\eps$-close rectangles. 

Let $R_Q$ be the set of all components of $R$ that contain a square. We may assume that all these components are circles, otherwise a component would come arbitrary close to the boundary of $P_4^\circ$, so there would be an $\eps$-close rectangle on it with aspect ratio $r$. If we could do this for all $\eps$, then a limit argument would give us a proper rectangle of aspect ratio $r$. So if need be, we choose a smaller $\eps$ for which this does not happen.

$R$ is a one-dimensional $\ZZ/4$-manifold, so $\ZZ/4$ acts on $R_Q$ as well. We decompose $R_Q=R_1\dot\cup R_2\dot\cup R_4$, where $R_1$ is the set of components with isotropy group $\langle0\rangle$, $R_2$ with isotropy group $\ZZ/2=\langle\eps^2\rangle\subset\ZZ/4$ and $R_4$ with $\ZZ/4$. Now we only need to count the number of squares on each $R_i$.

\begin{itemize}
\item $\sharp Q=4\mod 8$, since modulo $\ZZ/4$ it is odd (see Section \ref{secShnirelmansProof}). 
\item Every component $C\in R_Q$ contains an even number of squares, since while passing a square the rectangle changes from fat to skinny or vice versa (this follows from the bijectivity of the differential $\dd f$ at points in~$Q$).
\item 4 divides $\sharp R_1$, and every component in $R_1$ contains an even number of squares.
So the number of squares in the components of $R_1$ is divisible by~8.
\item 2 divides $\sharp R_2$, and if a component in $R_2$ contains a square~$S$, then it inscribes also $\eps^2\cdot S$.
When it goes through a square and changes from fat to skinny, then so it does at $\eps^2\cdot S$.
Hence it has to go through $4k$ squares, $k\geq 1$.
Thus the total number of squares in the components of $R_2$ is divisible by~8.
\item If a component $C$ of $R_4$ goes through a square $S$ and changes from fat to skinny, then it also goes through $\eps\cdot S$ and changes from skinny to fat. That is, in between it had to go through an even number of squares, all of which of course belong to a different $\ZZ/4$-orbit. Hence the number of square-orbits on $C$ is odd, $\sharp(Q\cap C)=4\mod 8$. 
\end{itemize}

Putting this modulo $8$ together, we get $\sharp R_4=1\mod 2$, which is even a bit stronger than what is stated in the lemma.
\end{proof}

\begin{remark}
Another viewpoint for the proof of Lemma~\ref{lemRectanglesOnCurvesUsingRectangles} is the following.
As a $1$-manifold, $R$ is a union of circles and segments.
The endpoint of the segments correspond to degenerate rectangles.
These rectangles are either ``skinny'' or ``fat'', and the segment $s$ has either endpoints of the same type or of different type; we say that $s$ type-preserving or type-reversing, respectively.
This feature of $s$ stays invariant under the action of~$\ZZ/4$.
One property is that type-preserving $s$ contains an even number of squares, and type-reversing $s$ an odd number.

By a bordism argument one can show that the number of type-reversing $s$ mod $\ZZ/4$ equals $1$ plus the number of $\ZZ/4$-invariant circle components of~$R$ modulo~$2$.
Moreover, type-reversing segments in $R$ contain rectangles of any aspect ratio.
Thus this yields another proof of the lemma.
\end{remark}

\subsection{Inscribed rectangles with aspect ratio $\sqrt{3}$} \label{secRectanglesProofOfMainTheorem}

In the case $r=\sqrt{3}$ there is a \qm{hidden} symmetry that we will use to prove Theorem~\ref{thmRectangularPegForSqrt3}.

\begin{figure}[htb]
  \centering
  \begin{minipage}[b]{0.35\textwidth}
	\centering
    \input{Figures/AngularConvexity2.pspdftex}
    \caption{Example of a curve that is not $60^\circ$-angular convex.}
    \label{figAngularConvexity}
  \end{minipage}
\quad
  \begin{minipage}[b]{0.55\textwidth}
    \centering
    \input{Figures/RectangularPegStern2.pspdftex}
    \caption{A star containing three $60^\circ$-parallelo\-grams; $x_1z_1x_2z_2$ is skinny, the other two are fat.}
    \label{figRectPegStern}
  \end{minipage}
\end{figure}

We leave all technical details concerning transversality to the subsequent Section~\ref{secTechnicalities}.
Suppose we are given a smooth curve $\gamma:S^1\incl \RR^2$. 
Define a map
\[
\begin{aligned}
f:\ &(S^1)^4&\To&_{G}&\ \RR^2\times S^1 \\
&(x_1,x_2,y_1,y_2)&\longmapsto&&\ (v,\alpha),
\end{aligned}
\]
where $v$ is again the difference of the diagonal midpoints in the quadrilateral $(\gamma(x_1),\gamma(y_1),\gamma(x_2),\gamma(y_2))$ and $\alpha$ is the mod-$180^\circ$ angle between these diagonals (we measure angles always in counter-clockwise sense). 
If one diagonal is degenerate to a point we take the tangent of $\gamma$ at this point to define $\alpha$.

The map $f$ is $G$-invariant, where $G:=\ZZ/2\times\ZZ/2=\{\bar{0}_x,\bar{1}_x\} \times \{\bar{0}_y,\bar{1}_y\}$ acts on $(S^1)^4$ by $\bar{1}_x\cdot(x_1,x_2,y_1,y_2)=(x_2,x_1,y_1,y_2)$ and $\bar{1}_y\cdot(x_1,x_2,y_1,y_2)=(x_1,x_2,y_2,y_1)$.

Let $P:=f^{-1}(0,60^\circ)$ be the set of inscribed parallelograms in $\gamma$ having a $60^\circ$-angle modulo $180^\circ$ between their diagonals. 
We call them $60^\circ$-parallelograms.
We may assume that $P$ is a union of connected $1$-dimensional submanifolds $K_i$ of $(S^1)^4$,
\[
P = K_1\cup\ldots\cup K_n,\ K_i\homeo S^1.
\]
$P$ does not contain points $(x_1,x_2,y_1,y_2)\in P$ with $x_1=x_2$ or $y_1=y_2$, since $\gamma$ was assumed to be $60^\circ$-angular convex. 
Thus, $G$ acts freely on $P$.
We denote $(S^1)^4/_G=M^2$ where $M:=(S^1)^2/_{\ZZ/2}$ is the M\"{o}bius strip. 
The first factor $M$ parametrises $x_1$ and $x_2$ without their order and the second $M$ parametrises $y_1$ and $y_2$. 
Let $L_1\cup\ldots\cup L_m\subset M^2$ be the quotient manifold $(\bigcup K_i)/_G$.
Then $L$ represents an element in the $1$-dimensional unoriented bordism group $\unorientedBordism_1(M^2)\iso \unorientedBordism_1((S^1)^2)\iso (\ZZ/2)^2$, since all $60^\circ$-angular convex curves are isotopic in the plane and $G$-homotopies of $f$ change $K_1\cup\ldots\cup K_n$ by a $G$-bordism.

If $\gamma$ is the unit circle then we see that $P$ is the disjoint union of two circles that get identified by~$G$. Their quotient $L$ is one circle that represents $(\bar{1},\bar{1})\in \unorientedBordism_1(M^2)\iso(\ZZ/2)^2$, where $\bar{1}\in\ZZ/2$ is the generator.

$P$ does not contain parallelograms that have an edge that is degenerate to a point. 
Hence the $x_1$ and $x_2$-coordinates will always differ from the $y_1$ and $y_2$-coordinates at any point $(x_1,x_2,y_1,y_2)\in P$. 
Therefore the circles $L_i$ can only represent the elements $(\bar{0},\bar{0})$ and $(\bar{1},\bar{1})$ of $\unorientedBordism_1(M^2)\iso (\ZZ/2)^2$.

Now we come to the \qm{hidden symmetry}, that is, the geometric piece of information that is the key in this proof. 
Let $W:=\{(\alpha,\beta,\gamma)\in (S^1)^3\st \alpha+\beta+\gamma=0^\circ\mod 180^\circ\}$.
We define a map
\[
\begin{aligned}
F:\ &(S^1)^6&\To&_{}&\ (\RR^2)^3\times W \\
&(x_1,x_2,y_1,y_2,z_1,z_2)&\longmapsto&&\ (m_x,m_y,m_z,\alpha_{yz},\alpha_{zx},\alpha_{xy}),
\end{aligned}
\]
where $m_x$ is the mid-point of the segment $(\gamma(x_1),\gamma(x_2))$, $\alpha_{yz}$ is the mod-$180^\circ$-angle between the lines through $(\gamma(y_1),\gamma(y_2))$ and $(\gamma(z_1),\gamma(z_2))$ (if $y_1=y_2$, we take instead the tangent of $\gamma$ at this point, and analogously in case $z_1=z_2$), and analogously for the the other coordinates. 
$F$~is equivariant with respect to the natural actions of the wreath product $G':=(\ZZ/2)^3\rtimes \ZZ_3$.
Let 
\[
\widetilde{S}:=F^{-1}(\diag_{(\RR^2)^3}\times \{(60^\circ,60^\circ,60^\circ)\}).
\]
We may assume that $\widetilde{S}$ is a $0$-dimensional free $G'$-manifold.
We call $S:=\widetilde{S}/_{G'}$ the set of \emph{stars}. 
Every star $s\in S$ contains three $60^\circ$-parallelograms on $\gamma$, namely $P_{yz}$, $P_{zx}$ and~$P_{xz}$, see Figure~\ref{figRectPegStern}. 
Modulo $G$ they lie in some components $L_i$, $L_j$, and~$L_k$ (they are not necessarily pairwise distinct).
We say that this star $s$ relates $L_i$, $L_j$, and~$L_k$. 
Saying this is unique up to cyclic permutation of $L_i$, $L_j$, and~$L_k$.
So we can draw a directed graph $D$ (with possibly multiple arcs and loops) whose nodes are the components of~$L$, and we draw for each star a directed triangle $L_i\to L_j\to L_k\to L_i$.

Assume that $\gamma$ does not inscribe a rectangle of aspect ratio~$\sqrt{3}$.
These are exactly the rectangles whose diagonals cross in a $60^\circ$-angle.
Then all $60^\circ$-parallelograms on $\gamma$ are \emph{skinny} or \emph{fat} in the sense that the $x$-diagonal is longer or shorter than the $y$-diagonal, respectively. 
By continuity this does not change along the components of~$L$.
Hence we can call the $L_i$'s fat or skinny.

Along a component $K_i$, $\{x_1,x_2\}$ and $\{y_1,y_2\}$ never intersect: If they did, by $m_x=m_y$ and $\alpha_{xy}=60^\circ$, all four points would need to coincide, but at the tangent at that point of $\gamma$ would also need to have two directions, which differ by $60^\circ$, which is impossible.
Thus for each $i$, $[L_i]\in \unorientedBordism_1(M^2)$ is either $(\bar{0},\bar{0})$ or~$(\bar{1},\bar{1})$. 
Correspondingly, we say that the \emph{winding number} $w(L_i)$ of $L_i$ is $\bar{0}$ (even) or $\bar{1}$ (odd), respectively.

Let $x,y: M^2\to M$ be the projections to the first and to the second factor, respectively.
An arc $L_i\to L_j$ in the graph $D$ corresponds to an intersection of $y(L_i)$ and $x(L_j)$.
The number of such intersections is 
\begin{equation}
\label{eqIntersectionNumberInMoebiusStrip}
\sharp(y(L_i)\cap x(L_j)) = w(L_i)\cdot w(L_j) \mod 2.
\end{equation}

We will derive a contradiction by double counting the number of stars $\sharp S$.

By \eqref{eqIntersectionNumberInMoebiusStrip}, components of $L$ with even winding number will have no influence on what follows.
Let $s$ be the number of skinny components of $L$ with odd winding number,
and let $f$ be the number of fat components of $L$ with odd winding number.

We know that $[L]=\sum_i [L_i]=(\bar{1},\bar{1})$, thus 
\[
s+f = 1 \mod 2.
\]
Note that no star relates three skinny or three fat $60^\circ$-parallelograms with each other.
Hence every star gives exactly one arc from a skinny to a fat component of~$L$.
Modulo~$2$ and using~\eqref{eqIntersectionNumberInMoebiusStrip}, there are congruent $s\cdot f=0 \mod 2$ of these arcs.
Therefore,
\[
\sharp S = 0 \mod 2.
\]

On the other hand, every star relates three components, two of which are skinny or two of which are fat.
So every star gives exactly one arc between two skinny components or between two fat components.
Using \eqref{eqIntersectionNumberInMoebiusStrip}, the number of arcs between skinny components modulo two is
\[
s^2 = s \mod 2,
\]
and the number of arcs between fat components modulo two is
\[
f^2 = f \mod 2.
\]
Together this gives,
\[
\sharp S = s+f = 1 \mod 2.
\]
This is a contradiction, which finishes the proof of Theorem~\ref{thmRectangularPegForSqrt3}.
\qed

\subsection{Technical Details} \label{secTechnicalities}

In the previous section we assumed that the set of inscribed $60^\circ$-parallelograms $P$ is a $1$-dimensional manifold in the $4$-manifold $M^2$.
Also the set of stars should be finite.
At the same time, when two parallelograms $p_1$ and $p_2$ have a common diagonal $y(p_1)=x(p_2)$ they form a star. Thus there should be another parallelogram $p_3$ such that $x(p_1)=y(p_3)$ and $y(p_2)=x(p_3)$. 
These triple intersection points come from the geometry, but they are in some sense \emph{not generic}.
That is, we need to be careful on how to make the test maps $f$ and $F$ simultaneously transversal in order to keep the geometric property of a star and without violating the equivariance.
We solve this issue by perturbing the following two maps.

Let 
\[
m: (S^1)^2\to \RR^2
\]
be the map that sends $(x_1,x_2)\in (S^1)^2$ to the mid-point $\frac{\gamma(x_1)+\gamma(x_2)}{2}$.
Let
\[
\alpha: (S^1)^2\to S^1
\]
be the map that sends $(x_1,x_2)\in (S^1)^2$ to the mod-$180^\circ$ angle of the line through $\gamma(x_1)$ and $\gamma(x_2)$ and some fixed line in the plane. The maps $f$ and $F$ can written in terms of $m$ and $\alpha$,
\[
f(x_1,x_2,y_1,y_2)=\left(m(y_1,y_2)-m(x_1,x_2),\alpha(y_1,y_2)-\alpha(x_1,x_2)\right)
\]
and similarly $F$.

Let $\varphi_i:S^1\to [0,1]$, $i=1\ldots k$, be a partition of unity of $S^1$ subordinate to a covering of $S^1$ with small $\eps$-balls.
We will perturb the maps $m$ and $\alpha$ with two sets of parameters 
$S_m:=([-\eps,+\eps]^2)^{\binom{k+1}{2}}$ 
and $S_\alpha:=[-\eps,+\eps]^{\binom{k+1}{2}}$ as follows:
\[
\begin{aligned}
m':\ &S_m\times (S^1)^2&\To&\ \RR^2 \\
&(s_m,x_1,x_2)&\longmapsto&\ m(x_1,x_2)+\sum_{i\leq j} (\varphi_i(x_1)\varphi_j(x_2)+\varphi_i(x_2)\varphi_j(x_1))\cdot(s_m)_{i,j},
\end{aligned}
\]
and 
\[
\begin{aligned}
\alpha':\ &S_\alpha\times (S^1)^2&\To&\ S^1 \\
&(s_\alpha,x_1,x_2)&\longmapsto&\ \alpha(x_1,x_2)+\sum_{i\leq j} (\varphi_i(x_1)\varphi_j(x_2)+\varphi_i(x_2)\varphi_j(x_1))\cdot(s_\alpha)_{i,j},
\end{aligned}
\]
This defined analogous functions $f':S_m\times S_\alpha\times (S^1)^4\toG{G} \RR^2\times S^1$ and $F':S_m\times S_\alpha\times (S^1)^6\toG{K} (\RR^2)^3\times W$.
Because of the additional parameter space $f'$ and $F'$ are transversal to the respective test-spaces $\{(0,60^\circ)\}$ and $\Delta_{(\RR^2)^3}\times\{60^\circ,60^\circ,60^\circ\}$.
By the transversality theorem \cite[p. 68]{GuPo10diffTop}, for almost all choices $s:=(s_m,s_\alpha)$ (up to a zero set), the perturbations $f'_s:=f'(s,\underline{\ \ })$ and $F'_s:=F'(s,\underline{\ \ })$ are transversal to the test-spaces as well. Similarly one can show that for almost all $s$, $y(K_i)$ intersects $x(K_j)$ transversally for all $i$, $j$.

\section{Inscribed crosspolytopes} \label{secCrosspolytopesOnSpheres}

In Klee \& Wagon \cite[Problem 11.5]{KlWa96problemsInPlaneGeomAndNumberTh} is was asked whether every $3$-dimensional convex body circumscribes the vertices of a regular octahedron.
Makeev \cite{Mak03universallyInscrAndOutscrPolytopes} proved this for smooth convex bodies and Karasev \cite{Kar09inscrRegularCrosspolytope} generalised the proof to smoothly embedded spheres in higher dimensions as follows.


\begin{theorem}[Makeev, Karasev]
Let $d$ be an odd prime power. 
Then every smooth embedding $\Gamma:S^{d-1}\to\RR^d$ contains the vertices of a regular $d$-dimensional cross\-polytope.
\end{theorem}

In 1965 H. Guggenheimer \cite{Gug65finiteSetsOnCurvesAndSurfaces} gave already a proof for all~$d$, however there is unfortunately an error in his main lemma due to some connectivity arguments, which seems to invalidate the proof. Recently, Akopyan and Karasev \cite{AkKa11inscrRegOctahedron} proved by a careful and non-trivial approximation argument that for $d=3$, the smooth embedding $\Gamma$ can be replaced by the boundary of a simple polytope.

An interesting possible extension seems to be the following.
Let $M$ be a Riemann manifold.
By an \textdef{inscribed crosspolytope} $P$ in $M$ we mean a set of $2d$ pairwise distinct points $v^{\eps}_{i}\in M$, $\eps\in\{+,-\}$, $i\in\{1,\ldots,d\}$.
We call two vertices $v^\eps_i$ and $v^\delta_j$ \textdef{opposite} if $i=j$ and $\eps=-\delta$.
If any pair of non-opposite vertices of $P$ have the same distance in $M$ then we call $P$ a \textdef{regular crosspolytope}.

\begin{conjecture}[``Crosspolytopal peg problem for manifolds'']
\label{conjCrosspolytopesOnManifolds}
Let $d$ be a positive integer.
Then every smooth embedding $\Gamma:S^{d-1}\to M$ into a Riemann manifold contains the vertices of a regular $d$-dimensional cross\-polytope.
\end{conjecture}

The aim of this section is to show that the conjecture in general is probably very difficult.

\newcommand{\RRtilda}{\mbox{\begin{minipage}[t]{0.8em}$\RR$\\ \vspace{-1.95em} $\sim$\end{minipage}}}

\paragraph{The topological counter-example.}
A solution of the conjecture would involve deeper geometric reasoning, since there is the following ``topological counter-example'' for $d=3$.
Suppose we are given a smooth embedding $\Gamma:S^2\to M$.
Let $G\iso(\ZZ/2)^3\rtimes S_3$ be the symmetry group of the regular octahedron and $G_{\rm{or}}\subset G$ be the subgroup of orientation preserving symmetries.
$G$ acts on $(S^2)^6$ by permuting the coordinates in the same way as it permutes the vertices of the regular octahedron.
Let $G$ act on $\RR^{12}$ by permuting the coordinates in the same way as it permutes the edges of the regular octahedron.
The subrepresentation $(\diag_{\RR^{12}})^\bot\subset\RR^{12}$ is denoted by $Y$.
Let $\fatdiag_{(S^2)^6}$ be the space of all $6$-tuples in $(S^2)^6$ that contain at least two equal elements, that is, the fat diagonal.
Let $B$ be a small $\eps$-neighborhood of $\fatdiag_{(S^2)^6}$, where $\eps$ depends only on an isotopy of $\Gamma$ to some nice embedding, that we will describe later.
Then the complement $X:=(S^2)^6\wo B$ is a free compact $G$-manifold with boundary and 
\[
X \homeq_G \{(x_1,\ldots,x_6)\in(S^2)^6\st x_i\textnormal{ are pairwise distinct}\}=(S^2)^6\wo\fatdiag_{(S^2)^6}.
\]
Then $\Gamma$ induces a test map
\[
t:X\toG{G} Y,
\]
which measures the lengths of the edges of the parametrised octahedra modulo $\one=(1,\ldots,1)$.
This map depends only on the distance function $d:M\times M \to \RR$ on $M$.
Since $\eps$ was chosen to be small, $t|_{\bd X}$ maps into $Y\wo\{0\}$ and this map is unique up to $G$-homotopy. 
We will use this fact later to assume that $\Gamma$ is actually some nice embedding of $S^{d-1}$ into $\RR^d$.
The solution set $S$ of regular octahedra inscribed in $\Gamma$, in the sense that all edge lengths coincide, is $S:=t^{-1}(0)$.
The subset $S_{\rm{or}}\subset S$ of positively oriented inscribed octahedra is a part of the preimage $t^{-1}(0)$, and $t$ induces an isomorphism of $G_{\rm{or}}$-vector bundles over $S_{\rm{or}}$,
\[
TS_{\rm{or}} \oplus (i_{S_{\rm{or}}})^*(X\times Y) \iso (i_{S_{\rm{or}}})^*(TX),
\]
where $i_{S_{\rm{or}}}$ denotes the inclusion $S_{\rm{or}}\incl X$, and $X\times Y$ is considered as the trivial vector bundle over $X$ with fiber~$Y$.
Thus $S_{\rm{or}}$ together with this normal data represents an element $[S_{\rm{or}}]$ in the equivariant normal bordism group (see Koschorke \cite[Chap. 2]{Kos81vectorFieldsAndOtherVectorBundleMorph})
\[
\Omega_1^{G_{\rm{or}}}(X,X\times Y-TX)=\Omega_1(X/_{G_{\rm{or}}},X\times_{G_{\rm{or}}} Y-T(X/_{G_{\rm{or}}})),
\]
which is well-defined, since $\ZZ/2$-homotopies of $d$ relative to a small neighborhood of $\Delta_{M^2}$ change $S$ only by a normal bordism that stays away from the $\bd X$ if $\eps$ was chosen small enough, and components of octahedra of different orientation are always separated from each other.
In Koschorke's notation, $[S_{\rm{or}}]$ is the obstruction 
\[
\widetilde{\omega}_1(\RRtilda,X\times_{G_{\rm{or}}} Y,(\id_{\bd X},t|_{\bd X})/_{G_{\rm{or}}}), 
\]
where $\RRtilda$ denotes the trivial line bundle.
\begin{theorem}
The above defined $[S_{\rm{or}}]$ is zero.
Hence 
\[
[S]\in\Omega_1^G(X,X\times Y-TX)
\]
is zero as well.
In particular, the test map t can be deformed $G$-equivariantly relative to $\bd X$ to a map $t'$, such that $0\not\in t'(X)$.
\end{theorem}
The existence of the test map $t'$ that fulfills the boundary conditions is what we call a topological counter-example.

\begin{proof}[Sketch of Proof]
To construct a convenient representative for $[S_{\rm{or}}]$ we take for $\Gamma$ a parametrisation of the  ellipsoid $\{(x,y,z)\st x^2+y^2+2z^2=8\}$, which has a rotational symmetry about the $z$-axis.
We let $t$ and $S$ be the corresponding test map and solution set, respectively.
We compute $S$ explicitly as a real algebraic variety; the relevant SageMath script is available at~\cite{Mat09squarePegRepository}. 
One inscribed regular octahedron has the six vertices $\pm(2,0,\sqrt{2})$ and $\pm(1,\pm\sqrt{3},-\sqrt{2})$.
If we rotate this octahedron around the $z$-axis then we get up to symmetry all inscribed octahedra in $\Gamma$, and $S$ is a disjoint union of $16=\frac{1}{3}\cdot\sharp G$ circles; see Figure~\ref{figInscribedOctahedron}.
\begin{figure}[h]
\centering 
\input{Figures/SpecialOctahedronOnSphere.pspdftex}
\caption{Orthogonal projection of $\Gamma$ into the $xy$-plane together with an inscribed regular octahedron.}
\label{figInscribedOctahedron}
\end{figure}

The $G$-bundles $X\times Y$ and $TX$ are $G$-orientable.
Therefore the relevant part of Koschorke's exact sequence \cite[Thm. 9.3]{Kos81vectorFieldsAndOtherVectorBundleMorph} becomes
\[
\begin{aligned}
H_2(X/_{G_{\rm{or}}};\ZZ)\to \ZZ/2\to \Omega_1(X/_{G_{\rm{or}}},& X\times_{G_{\rm{or}}}Y-T(X/_{G_{\rm{or}}}))
\\ &
\to H_1(X/_{G_{\rm{or}}};\ZZ)\to 0.
\end{aligned}
\]
It is not difficult to see that the image of $[S_{\rm{or}}]$ in $H_1(X/_{G_{\rm{or}}};\ZZ)=H_1(G_{\rm{or}};\ZZ)$ is zero.
This is because the 120 degrees rotation of a regular octahedron about the line connecting the midpoints of two opposite triangles is an element of the commutator of $G_{\rm{or}}$.
It requires more visualisation to see that $[S_{\rm{or}}]$ is in fact the image of the generator of $\ZZ/2$.
The hard part is to show that $\ZZ/2$ unfortunately lies in the image of $H_2(X/_{G_{\rm{or}}};\ZZ)$, which we computed with a rather long computer program, see~\cite{Mat09squarePegRepository}. %
It finds that $H_2(X/_{G_{\rm{or}}};\ZZ)\iso \ZZ/4\times (\ZZ/2)^3$, where one can choose the generators such that the first three map to zero and the last one to the generator of $\ZZ/2$.

The $G_{\rm{or}}$-null-bordism of $S_{\rm{or}}$ can be extended to a $G$-null-bordism of $S$.
By Theorem 3.1 of Koschorke \cite{Kos81vectorFieldsAndOtherVectorBundleMorph}, we can extend the section as stated.
\end{proof}

\paragraph{Remarks to the algorithm.} 
An economical $S_6$-CW-complex structure on $(S^2)^6$ is based on an $S_6$-cell decomposition of $\RR^2$ of Fuks~\cite{Fuks70cohomologyOfBraidGroupsMod2} and Vassiliev~\cite{Vas88braidGroupCohomAndAlgComplexity}, which has few high dimensional cells.
$\fatdiag_{(S^2)^6}$ is a subcomplex, so one can compute $$H_2(X/_{G_{\rm{or}}})\iso H^{10}((S^2)^6/_{G_{\rm{or}}},(\fatdiag_{(S^2)^6})/_{G_{\rm{or}}}).$$
The Smith normal form is used to compute this cellular cohomology and the LLL-algorithm to choose economical generators.
The image in $\ZZ/2$ is determined by computing second Stiefel-Whitney classes, which have been implemented as obstruction classes.

\section*{Acknowledgements}

I want to thank to Karim Adiprasito, Ulrich Bauer, Pavle \blagojevic, Jason Cantarella, Roman Karasev, Sergey Melikhov, Carsten Schultz, John Sullivan, Helge Tverberg, Sini\v{s}a \vrecica, G\"unter Ziegler, Aleksey Zinger, and Rade \zivaljevic\ for many very useful discussions. Moreover I thank the referees for very valuable comments and questions.

This work was supported by Studienstiftung des dt.~Volkes and Deutsche Telekom Stiftung at Technische and at Freie Universit\"at Berlin, NSF Grant DMS-0635607 at Institute for Advanced Study, by an EPDI fellowship at Institut des Hautes Etudes Scientifiques, Forschungsinstitut f\"ur Mathematik, and Isaac Newton Institute, by Max Planck Institute for Mathematics Bonn, and by Simons Foundation grant {\#}550023 at Boston University.

\small
\bibliographystyle{plain}

\bibliography{../mybib07}

\vspace{0.5cm}
\noindent
\small{
Boston University \\
matschke@bu.edu
}
\end{document}